\documentclass[12pt,a4paper]{article}
\usepackage{amssymb}
\usepackage{amsmath,amsfonts,amsthm,color}
\usepackage{graphicx}
\usepackage{amsmath}
\usepackage{amsthm,amscd}
\usepackage{graphicx}
\usepackage{amsmath}
\usepackage{amstext}
\usepackage{hyperref}

\newtheorem{maintheorem}{Theorem}

\newtheorem{theorem}{Theorem}
\newtheorem{corollary}[theorem]{Corollary}
\newtheorem{definition}[theorem]{Definition}
\newtheorem{example}[theorem]{Example}
\newtheorem{lemma}[theorem]{Lemma}

\newtheorem{proposition}[theorem]{Proposition}
\newtheorem{remark}[theorem]{Remark}

\title{Thermodynamic Formalism for Non-Invariant Probabilities: Entropy, Divergence, Channels}

\author{Artur O. Lopes and Elismar R. Oliveira\\IME - UFRGS - Porto Alegre - RS - Brazil}

\begin{document}

\maketitle

\begin{abstract}
Let $\mathcal{M}$ and $\mathcal{M}_\sigma$ denote, respectively, the sets of Borel and $\sigma$-invariant probabilities on the uncountable set $\Omega=\{1,2,...,d\}^\mathbb{N}$. Given a H\"older Jacobian $J$,  the dual Ruelle operator $\mu \in \mathcal{M} \mapsto \mathcal{L}_{J}^* (\mu)$ will play the role of a dynamical channel; generalizes the action of a stochastic matrix on probability vectors. Since $\mathcal{L}_{J}^*$ generally breaks $\sigma$-invariance, we introduce the \emph{ni-entropy} $h_{ni}(\mu)$ for non-invariant measures, employing the concept of the expected value function on cylinders of size one. $h_{ni}(\mu)$ differs from Shannon entropy but extends the Kolmogorov-Sinai entropy from $\mathcal{M}_\sigma$ to $\mathcal{M}$; $h_{ni}$ is concave and u.s.c. This leads us to the concept of ni-pressure $\mathcal{P}(A)$ of a potential $A:\Omega \to \mathbb{R}$. $\mu\in \mathcal{M}$ is not a density with respect to a dominating measure. We define an {\it ni-Kullback-Leibler divergence} $D_{ni}(\rho || \nu)$ not assuming 
absolute continuity. We investigate concepts of Information Theory such as Minimal Output Entropy, Donsker-Varadhan variational principle, a variational formula for $h_{ni}$, bounds for divergence, information projections, and  their connections to the ni-pressure $\mathcal{P}$.
\end{abstract}

\section{Introduction} \label{sec:Intro}

We introduce a natural concept of entropy $h_{ni}(\mu)$  for probabilities $\mu$ on the symbolic space $\Omega = \{1,2,\dots,d\}^\mathbb{N}$ that are not necessarily shift-invariant, where the concept of the expected value function plays an important role, along with an associated ni-variational pressure problem. Unlike Classical Information Theory, our framework accommodates classical channels acting on probabilities without assuming a dominating measure or absolute continuity; our setting is non-dynamical on the one hand, since we do not consider averages over finite orbits under the shift iteration, but somewhat dynamical on the other,  due to the use of the dual of the Ruelle operator as a channel (which generalizes  the action of a line stochastic matrix acting on the left on vectors of probabilities as mentioned in Remark 1 Chapter 2 page 27 in \cite{PP}).

\begin{definition}
	Let $\Omega=\{1,2,\dots,d\}^\mathbb{N}$ ($d\geq 2$) be the compact shift-space endowed with the metric $d(x,y)=2^{-k}$, where $k=\min\{i\geq 1 \mid x_i \neq y_i\}$, for $x \neq y$ and 0 otherwise. Let $\sigma: \Omega \to \Omega$ be the shift map $\sigma(x_1,x_2,\dots)= (x_2,x_3,\dots)$, and $\mathcal{F}$ the Borel $\sigma$-algebra. We denote by $\mathcal{M}$ the space of Borel probabilities on $\Omega$ (compact in the weak-$*$ topology), and by $\mathcal{M}_\sigma$ the subset of $\sigma$-invariant probabilities,  corresponding in Statistical  Physics to stationarity.
\end{definition}

$\{1,\dots,d\}$ is the space of symbols. For $a_1,\dots,a_k \in \{1,\dots,d\}$,  $\overline{a_1,\dots,a_k}= \{(a_1,\dots,a_k,x_{k+1},\dots) \mid x_j \in \{1,\dots,d\}, j\geq k+1\}$ is called a cylinder. 
Let $\mathcal{X}_1 = \sigma^{-1} (\mathcal{F})$. Given $j$ and $\rho\in \mathcal{M}$, denote $E_{\rho} (I_{\overline{j}}\mid \mathcal{X}_1)$  the conditional expectation of the indicator function $I_{\overline{j}}$ with respect to $\mathcal{X}_1$.  We use the convention $0\log 0=0$.
Our entropy is defined by \eqref{fifa1}, and we call it the ni-entropy of  $\rho\in \mathcal{M}$, distinguishing it from the classical meaning $h(\rho)$ for $\rho \in \mathcal{M}_\sigma$:
\begin{equation} \label{fifa1}
 h_{ni}(\rho) =	- \int \sum_{j=1}^d E_{\rho} (I_{\overline{j}} \mid \mathcal{X}_1) \log (E_{\rho} (I_{\overline{j}} \mid \mathcal{X}_1)) \, d\rho.
\end{equation}

Kolmogorov-Sinai entropy $h(\rho)$ can also be described by \eqref{fifa1} (see Section 3 in \cite{PP}, or the discussion on page 120 in \cite{Mart}, which describes the relevance of the concept of entropy defined in such a way in Information Theory).

Viewing our setting from the perspective of a Stochastic Process $X_n$, $n=0, 1, 2, \dots$, taking values in $\{1,2,...,d\}$ - not necessarily stationary - let $\mu$ denote the associated probability measure on $\{1,2,...,d\}^\mathbb{N}$.  $E_{\mu} ( I_{\overline{j}} | \mathcal{X}_1))$ should be understood as the conditional probability of $X_0=j$, given the information $(X_1, X_2, \dots, X_n, \dots)$. If $\mu$ were a Markov measure, this would correspond to $\mu(X_0=j|X_1)$.  In this case, in Example \ref{gost1} the function $E_{\mu} ( I_{\overline{j}} | \mathcal{X}_1)$ was explicitly computed for $j=1,2$, and $h_{ni}(\mu)$ is the mean information.  For the stationary case, the  entropy rate (6.4.1) in Sec. 6.3 \cite{Moser} is  aligned with \eqref{fifa1}.  

\begin{remark}Our results lead to expressions that are analogous to, but not exactly the same as, the classical case (which assumes a dominating measure). The upcoming expressions \eqref{cum}, \eqref{cum1}, \eqref{cum3},  dealing with the cumulant function and global bounds for divergence, and \eqref{triuvi87}, dealing with the tilting probability associated with information projections,  illustrate these differences.

\end{remark}

\begin{definition}
	Let $\mathcal{C}$, $\mathcal{C}^+$, and $\mathcal{H}$ denote the spaces of continuous, positive continuous, and H\"older continuous functions $f: \Omega \to \mathbb{R}$, respectively. A  function $J \in \mathcal{C}^+$ is a Jacobian if $\sum_i J(i*x)=1$ for all $x\in \Omega$, where $i*x = (i,x_1,x_2,\dots)$. We denote by $\mathfrak{J}$ the set of $J$ such that $\log J\in \mathcal{H}$.
\end{definition}

\begin{definition}
	The Ruelle operator $\mathcal{L}_J: \mathcal{C} \to \mathcal{C}$ associated with a Jacobian $J\in \mathfrak{J}$ is $\mathcal{L}_J(f)(x) = \sum_i J(i*x)f(i*x)$. Its dual $\mathcal{L}_J^*: \mathcal{M} \to \mathcal{M}$ satisfies $\int_{\Omega} f \, d\mathcal{L}_J^*(\rho) = \int_{\Omega} \mathcal{L}_J (f) \, d\rho$, for  all $f \in \mathcal{C}$. Our results will be just for $J\in \mathfrak{J}.$
\end{definition}

 In Statistical Mechanics, $\mathcal{L}_J$ acts as the transfer operator. In dynamics it is a crucial tool for studying the ergodic properties of equilibrium states \cite{PP}.  Note that $\mathcal{L}_J(1)=1$. The action of $f \to \mathcal{L}_J(f)$ is a generalization of the one  of a $d\times d$ row-stochastic matrix $P$ on column vectors $v\in \mathbb{R}^d$ on its right-hand side, that is $v \to  P v$. For a fixed $x \in \Omega$, the Jacobians $J(i*x) = p(i)=p_x(i)$ - where $i \in \{1, 2, \dots, d\}$ - can be viewed as a probability distribution over the preimages of $x$. From the definition $\mathcal{L}_J(f)(x) = \sum_i p_x(i) f(i*x)$.
On the other hand, the action of the dual of the Ruelle operator $\mathcal{L}_J^*$ generalizes the action of a  row-stochastic matrix $P$ on probability row vectors $q$ on its left-hand side, that is, $q \to Pq$.
Here, we view $\mathcal{L}_J^*$ as a classical channel mapping input information $\mu_1$ to output $\mu_2= \mathcal{L}_J^*(\mu_1)$. $\mathcal{L}_J^*$ is a $1$-Wasserstein contraction when acting on $\mathcal{M}$, with a unique fixed point $\mu_J \in \mathcal{M}_\sigma$ \cite{KLS}.  One can show that the Kolmogorov-Sinai entropy $h(\mu_J)=-\int \log J\, d \mu_J$ \cite{PP}.

For $A \in \mathcal{C}$, the topological pressure is
\begin{equation} \label{AA}
	P(A):= \sup_{\rho \in \mathcal{M}_\sigma} \left( \int A \, d\rho + h(\rho)\right),
\end{equation}
where $h(\rho)$ is the Kolmogorov-Sinai entropy. A probability $\mu_A \in \mathcal{M}_\sigma$ maximizing \eqref{AA} is called an equilibrium probability for $A$ (and is unique for $A \in \mathcal{H}$ \cite{PP}), which aligns with the Second Law of Thermodynamics via the so-called MaxEnt method. The set $\mathcal{G}$ of all such equilibrium probabilities $\mu_A$, for $A\in \mathcal{H}$, precisely coincides with the fixed points of $\mathcal{L}_J^*$ for all $J \in \mathfrak{J}$. Probabilities in $\mathcal{M}_\sigma$ can be approximated by probabilities in $\mathcal{G}$ weakly and also in entropy (see \cite{Lo2}), a phenomena called ergodic abundance (see Proposition 2.36 and  Corollary 2.39 in \cite{Bru}). The entropy as a thermodynamic limit, classical Kullback-Leibler divergence  and its relation with the dual of the Ruelle operator appear in Section 2.9.2 in \cite{Bru}.
Previous work relating  Thermodynamic Formalism and Information Theory appeared in \cite{Ermo,Verb,Bag}.

The push-forward map $\mathfrak{T}:\mathcal{M} \to \mathcal{M}$, defined by $\int f \, d\mathfrak{T}(\mu) = \int (f \circ \sigma) \, d\mu$, serves as a recovery channel for $\mathcal{L}_J^*$. Specifically (see \cite{LO1}), for any $n \in \mathbb{N}$:
\begin{equation} \label{mai7}
	\text{if } (\mathcal{L}_J^*)^n(\nu)=\mu, \text{ then } \mathfrak{T}^n(\mu)= \nu.
\end{equation}
This makes $\mathcal{L}_J^*$ a Sufficient Classical Channel and $\mathfrak{T}$ its Decoder \cite{Cs,Mac}. No information is lost by the action of $\mathcal{L}_J^*$, as the original probability can be perfectly reconstructed; this ties into Sufficient Statistics (see \cite{Galla}). Further references on classical Information Theory include  \cite{PoWu,Cover, Gray}.   Section 3 in \cite{Mart} analyzes the interplay between dynamical entropy and Information Theory.

In \cite{Fein} the authors show a remarkable interplay of
the Hahn-Banach Theorem and the Second Law of Thermodynamics. In Section 1.2 in \cite{Fein} the authors ask the question:

{\it Query: In classical thermodynamics, is entropy defined only on
equilibrium states? If so, why?}

Digressing on this issue, the authors follow an interesting mathematical rationale based on the fundamental laws of thermodynamics, specifically on the understanding of the Kelvin-Planck statement and the Clausius-Duhem inequality. Finally, in Remark 6.3 they conclude that:

{\it For existence of entropy and thermodynamic temperature, equilibrium is irrelevant.}

Regarding the above, it makes sense to investigate the meaning of entropy for not necessarily $\sigma$-invariant probabilities; for instance, eigen-probabilities for the dual of the Ruelle operator.
Our entropy \textbf{is not} of the standard static form $S(\varphi)=-\int \varphi \log \varphi \, d\mu_{\text{ref}}$. The displayed density-based representation requires a dominating  measure  $\mu_{\text{ref}}$  with respect to which the density exists.

\begin{remark} \label{Fra}
	Many standard proofs in Information Theory (e.g.,  Theorem 2.3.2 in \cite{Gray} regarding Kullback-Leibler divergence $D(p||q)= \int p(x) \log \frac{p(x)}{q(x))} dx$ and the Donsker-Varadhan variational principle) rely on mutually absolutely continuous probabilities and cannot be directly applied to our symbolic space setting, where probabilities may be mutually singular.
	A natural question is how classical concepts involving absolutely continuous probabilities $p  \ll q$, like $p(x) \log \frac{p(x)}{q(x))} dx$, correspond to their analogues in our setting. We address this issue in Remark \ref{Fra1}.
	Using the expected value function, we present the result corresponding to this variational principle in \eqref{kloe4} of Section \ref{sec:ni-entropy}. Avoiding a dominating measure drastically alters computations (see Examples \ref{reco},\,\ref{orqi8}).
\end{remark}

To our knowledge, our ni-entropy definition \eqref{fifa1} for not necessarily $\sigma$-invariant probabilities is novel. Unlike classical settings, we have to consider  the expected value function, which can change discontinuously (Remark \ref{hjq}), meaning $h_{ni}$ is not everywhere continuous; it is upper semi-continuous and concave \eqref{bufer}. Our non-negative entropy differs from the one described in \cite{BCMV}.

	The paper is organized as follows. At the beginning of Section \ref{ChaL} we review Classical Channels  acting on probabilities over $\{1,\dots,d\}$, as a preliminary step to studying the case  on $\{1,\dots,d\}^\mathbb{N}$. Section \ref{ChaL} summarizes the main results on ni-entropy, the dual Ruelle operator as a channel, and ni-pressure. Section \ref{sec:ni-entropy} formalizes ni-entropy, its transformation under the dual Ruelle operator, and also the ni-Kullback-Leibler divergence.  For information projections see Section \ref{inpro}. Section \ref{pre}  computes probabilities maximizing ni-pressure. Section \ref{sec:Chan Cap} analyzes the Minimal and Maximal Output entropies of the $\mathcal{L}_J^*$ channel.  For applications and explicit examples see Section \ref{sec:Examples} . Section \ref{Conc} summarizes our results and mentions topics for future investigation.

We thank Juliano Tresoldi for helpful conversations on related topics.

\section{The ni-pressure and the dual of the Ruelle operator seen as a Classical Channel} \label{ChaL}

We start with a motivation from  Classical Information Theory.
A classical channel $\Psi$ maps probabilities to probabilities. A standard memoryless example (see \cite{Gray,PoWu}) is the action of a row-stochastic matrix $P$ on a non-dynamical probability vector $p=(p_1,p_2)$ (where $p_1, p_2 > 0$ and $p_1+p_2=1$):
\begin{equation} \label{Sa1}
	\Psi(p) = \begin{pmatrix} p_1 & p_2 \end{pmatrix} P = \begin{pmatrix} p_1 & p_2 \end{pmatrix} \begin{pmatrix} P_{11} & P_{12}\\ P_{21} & P_{22} \end{pmatrix} = \begin{pmatrix} q_1 & q_2 \end{pmatrix}. 
\end{equation}
In this context, the entries $P_{ij}$ represent the probability of corruption when transmitting state $i$ to $j$ (see \cite{Feins1, PoWu}). The action  $\mu \to \mathcal{L}_J^*$ generalizes \eqref{Sa1}.

The Shannon (or static) entropy for this non-dynamical setting is:
\begin{equation} \label{Sa3} 
	S(p) = - \sum_{j=1}^d p_j \log p_j.
\end{equation}
Unlike the dynamical Kolmogorov-Sinai entropy, which is affine \cite{Walters}, the Shannon entropy is just concave.
A key question in non-dynamical settings is how the channel alters information, which requires comparing the input and output entropies:
$S(p)$ and  $S(\Psi(p))$.
When $\Psi$ is a Markov channel defined by a \emph{doubly} stochastic matrix $P$, one can show (see \cite{Sa1}) that:
\begin{equation} \label{Sa11} 
	S(p) \leq S(\Psi(p)).
\end{equation}
This monotonicity reflects a thermodynamic operation where the channel makes the distribution more random (see \cite{Sa1}). 
If $P$ is not doubly stochastic, inequality \eqref{Sa11} \textbf{may not be true} (for example, if $P_{11}=1/3$ and $P_{21}= 9/10$).

The dual of the Ruelle operator for a normalized potential acts as the infinite-dimensional counterpart of a row-stochastic matrix $P$ (see  Remark I of Chapter 2 in \cite{PP}); the only point where dynamics somehow manifests itself is via $\mathcal{L}_J^*$. We study probabilities on the symbolic space $\Omega=\{1,2,\dots,d\}^\mathbb{N}$.

We view the mapping $\mu \to \Phi(\mu) := \mathcal{L}_J^*(\mu)$ as a classical channel, where $J \in \mathfrak{J}$; dynamics play a role via the inverse branches of $\sigma$. The maximal entropy measure is not in the image of $\mathcal{L}_J^*$ (see Example \ref{mai49}) for most $J\in \mathfrak{J}$.

\begin{remark} \label{Su1}
Given $J\in \mathfrak{J}$, for	 $\mu \in \mathcal{M}$, the image $\mathcal{L}_J^*(\mu)$ is not a $\sigma$-invariant probability, even if $\mu \in \mathcal{M}_\sigma$ (see Example \ref{mai}), unless $\mu$ satisfies $\mathcal{L}_J^*(\mu)=\mu$.
\end{remark}

To study problems analogous to \eqref{Sa11}  on the full space $\mathcal{M}$, we choose to analyze the action of the channel $\mathcal{L}_J^*$ and  to consider the ni-entropy $h_{ni}(\mu)$ (introduced in \eqref{fifa1}), due to Remark \ref{Su1}. Unlike static entropy (using expressions like \eqref{Sa3}), the new  concept  $h_{ni}$ requires a conditional expectation integral.

The ni-entropy $h_{ni}$ is non-negative, upper semicontinuous, bounded by $\log d$, and coincides with  classical Kolmogorov-Sinai entropy when $\mu \in \mathcal{M}_\sigma$ (Proposition \ref{iimp}). Examples \ref{mojic} and \ref{maili7} illustrate how dynamics can determine, in some cases, whether $h_{ni}$ is positive or not, when considering non-invariant probabilities. We show that $h_{ni}(\delta_x) = 0$ (Example \ref{moj16}). Measures supported on finite sets can eventually attain the maximal entropy value (Example \ref{mojic}).

Denote by  $L_1: \mathcal{C}^+\to \mathcal{C}^+$ the operator s.t.  $L_1 (u)(x) =v(x)=  \sum_{j=1}^d u(j,x)$. 
By Theorem \ref{tirri} in Section \ref{sec:ni-entropy}, we can alternatively express $h_{ni}(\rho)$  as: 
\begin{maintheorem} \label{uu} Given $\rho \in \mathcal{M}$, then 
	\begin{equation} \label{este5}
	h_{ni}(\rho) = \inf_{u \in \mathcal{C}^+} \int \log \left( \frac{L_1 (u) \circ \sigma}{u} \right) d\rho= \inf_{g \in \mathcal{C}} \int \log \left( \frac{L_1 (e^g) \circ \sigma}{e^g} \right) d\rho.
	\end{equation}
\end{maintheorem}

The next claim broadens results from \cite{LMMS} to $\rho \in\mathcal{G}$ (see also \cite{Lo7} and \cite{LO7}). From \eqref{este5} one can show that $\rho \to h_{ni}(\rho)$ is concave and upper semicontinuous.
\begin{corollary} \label{uute} Given $\rho \in \mathcal{M}_\sigma$, then 
	\begin{equation} \label{este55}
	h(\rho) = \inf_{u \in \mathcal{C}^+} \int \log \left( \frac{L_1 (u)}{u} \right) d\rho= \inf_{g \in \mathcal{C}} \int \log \left( \frac{L_1 (e^g)}{e^g} \right) d\rho.
	\end{equation}
\end{corollary}

Entropy can be obtained from a Ruelle operator; using $\mathcal{L}_{1/d}$ instead of $L_1$ we get $h(\rho)-\log d\leq 0$  (results for a more general space of symbols \cite{LMMS}).

Subsection \ref{nini} computes $h_{ni}(\nu)$ for various $\nu \in \mathcal{M}$, such as non-stationary Markov Processes (Example \ref{gost1}) and  eigenprobabilities of $\mathcal{L}_B^*$, $B:\Omega \to \mathbb{R}$ (see definition in  Example \ref{eigp}). We also estimate the ni-entropy for tilting probabilities  (Lemma \ref{ms1} and Example \ref{orqi}), also known as canonical probabilities (Definition 1.10 in \cite{FriVe}), which are associated with Information Projections.

Theorem \ref{rreT} computes explicitly the channel output ni-entropy.
\begin{maintheorem} \label{rreT}
	Given $J\in \mathfrak{J}$ and $\rho \in \mathcal{M}$
	$$
	h_{ni}(\mathcal{L}_J^*(\rho)) = - \sum_{j=1}^d \int J (j, x_1,x_2,\dots) \log J (j, x_1,x_2,\dots)\, d \rho(x_1,x_2,\dots),
	$$
\end{maintheorem}

Thus, $h_{ni}(\mathcal{L}_J^*(\rho))$ is affine in $\rho$; Theorem \ref{rreT} is proved in Theorem \ref{iimp}.

For some $\mu_1\in \mathcal{M}$, we provide explicit expressions for $\mathcal{L}_J^*(\mu_1)$  in Subsection \ref{runi}.
The property corresponding  to \eqref{Sa11}  does not hold in our setting for $J \in \mathfrak{J}$, because  the inequality $h_{ni}(\mathcal{L}_J^*(\mu)) \geq h_{ni}(\mu)$ is not universally valid. $\mu \to \mathcal{L}_J^*(\mu)$ can either increase or decrease entropy (Examples \ref{mai49} and \ref{moj25}).

In Definition \ref{KL}  we introduce the ni-Kullback-Leibler divergence $D_{ni}(\rho || \nu)$, which notably does not assume absolute continuity between $\rho$ and $\nu$. Subsection \ref{KLD} provides explicit computations. Some dissimilarities between the Kullback-Leibler divergence and the ni-Kullback-Leibler divergence   are described in Remarks \ref{Fra1} and \ref{moj16} . 
Subsection \ref{JS}  introduces the ni-Jensen-Shannon divergence providing a pseudo-metric on $\mathcal{M}$, and a metric on $\mathcal{G}.$

In our framework it is natural to introduce two variational problems for a continuous potential $A: \Omega \to \mathbb{R}$. The first evaluates channel output entropy:
\begin{equation} \label{AABB} 
	\mathfrak{P}(A) := \sup_{\mu \in \mathcal{M}} \left( \int A \, d\mu + h_{ni}(\mathcal{L}_J^*\mu) \right).
\end{equation}
A measure attaining this supremum is called a $\mathfrak{P}(A)$-maximizing probability.

The second, to be investigated more closely here, is the ni-pressure of $A$:
\begin{equation} \label{AAB} 
	\mathcal{P}(A) := \sup_{\mu \in \mathcal{M}} \left( \int A \, d\mu + h_{ni}(\mu) \right).
\end{equation}

A probability $\mu \in \mathcal{M}$ maximizing \eqref{AAB} is called a $\mathcal{P}(A)$-maximizing probability. The upper semicontinuity of $\mu \to h_{ni}(\mu)$ guarantees such a measure exists. Since Proposition \ref{iimp} ensures $h_{ni}(\mu)\leq \log d$ for all $\mu \in \mathcal{M}$, the maximal entropy measure is a $\mathcal{P}(0)$-maximizing probability, though not necessarily unique (see Example \ref{mojic}). For   H\"older potentials $A$, Examples \ref{foie} and \ref{foie1} show cases where no probability in $\mathcal{M}_\sigma$ maximizes $\mathcal{P}(A)$.
In Section \ref{inpro} we analyze properties of Information Projections (an overview in \cite{Nie}) obtained from the ni-Kullback-Leibler divergence and its relation to the ni-pressure $\mathcal{P}$.




The function $A \to \mathcal{P}(A)$ satisfies the abstract properties of Pressure outlined in Definition 2.2 of \cite{BCMV}, even though the entropy considered in \cite{BCMV} differs from ours.  Furthermore, adapting the proof of Theorem 9.12 in \cite{Walters},  for $\mu \in \mathcal{M}$ one can easily show that $h_{ni}(\mu)$ is the Legendre transform of $\mathcal{P}$:
\begin{equation} \label{B33} 
h_{ni}(\mu)=\inf_{g \in \mathcal{C}} \left( \mathcal{P} (g) - \int g d \mu   \right).
\end{equation}

We also analyze another variational problem: for a continuous $A:\Omega \to \mathbb{R}$:
\begin{equation} \label{AAB33} 
	\mathbb{P}(A)=\sup_{\mu \in \mathcal{M}} \left( \int A \,d \mathfrak{T}(\mu) + h_{ni}(\mu) \right).
\end{equation}

The inequalities $\mathcal P(A)\ge P(A)$ and $\mathfrak P(A)\ge P(A)$ follow from restricting/expanding the optimization domain only after using $h_{ni}=h$ on $\mathcal M_\sigma$. 
One can show $\mathbb{P}(A) \geq \mathcal{P}(A)$ (see Theorem \ref{tos}).

The thermodynamic relation between free energy (minus pressure) and mutual information is discussed in \cite{Feins1} (page 42). In this vein, Theorem \ref{bela} links the ni-pressure $\mathfrak{P}$ with concepts from Information Theory.

We prove in Theorem \ref{PreA21q} the following result for $J\in \mathfrak{J}$:
\begin{maintheorem} \label{PreA21}
	Take $A(x)=\log J(x_1,x_2,\dots)$. Then, $\mathcal{P}(A)=$
	$$ \sup_{\mu \in \mathcal{M}} \left( \int \sum_{j=1}^d I_{\overline{j}}(x_1,x_2,\dots)\log J (j,x_2,\dots,x_n,\dots)\,d \mu(x) + h_{ni} (\mu )\right)=0.$$
	The probability $\mu_J\in \mathcal{G}$ attains the value $\mathcal{P}(A)$; there are others not in $\mathcal{M}_\sigma$.
\end{maintheorem}
To address the variational problem $\mathfrak{P}(A)$, analyzing the range of the channel $\Phi=\mathcal{L}_J^*$ is crucial. In Section \ref{sec:Chan Cap}, among other results, we investigate finding $\rho$ to optimize various output entropies for a given $J \in \mathfrak{J}$ such as:
\begin{align} 
	H_i(\Phi) &= \inf_{\rho\in \mathcal{M}} h_{ni}( \mathcal{L}_{J}^* (\rho))  &\text{(Minimal output entropy)} \label{al1}\\
	H_s(\Phi) &= \sup_{\rho\in \mathcal{M}} h_{ni}( \mathcal{L}_{J}^* (\rho)) \quad &\text{(Maximal output of a channel)} \label{al2}
\end{align}

Problem $H_s(\Phi)$ is addressed in  Corollary 5.10 on page 72 in \cite{PoWu}. 
Theorem \ref{rreT} distinguishes clearly between optimization over $\mathcal M$ and over $\mathcal M_\sigma$, and identifies the latter with an ergodic-optimization problem.

We also analyze a constrained problem for a H\"older potential $A$. Fixing $k\in \mathbb{R}$, we set $\mathfrak{K}(k)= \{\mu \in \mathcal{M} \mid \int A \, d\mu=k\}$ and consider $H(\Phi,A,k) =$
\begin{equation} \label{pluu25} 
	 \sup_{\rho\in \mathfrak{K}(k) } h_{ni}( \mathcal{L}_{J}^* (\rho)) \quad \text{(a dynamical capacity cost problem)}.
\end{equation}

Problem \eqref{pluu25} considers questions addressed in Example 5.2 in \cite{PoWu}.
\begin{maintheorem} \label{bela} 
	Given the ni-pressure problem $\mathfrak{P}(A)$, there exists $k$ such that the probability maximizing $H_s(\mathcal{L}_{J}^*,A,k)$ also maximizes $\mathfrak{P}(A)$.
\end{maintheorem}
This will be proved in Theorem \ref{tos} and it is related to \cite{CL}.
\smallskip

Denote $\Lambda=\{q=(q_1,q_2,...,q_d), \sum_{j=1}^d q_j=1, q_j\geq 0, j=1,...,d\}$.
In  Subsection \ref{anot} we follow a different point of view: in some sense the channel $\mathcal{L}_J^* $ will act on probabilities $q\in \Lambda$ and not on $\mathcal{M}.$
Suppose $\mathcal{L}_J^* (\mu)=\mu$. By definition, the mutual information of $\mathcal{L}_J^*$ acting  on a probability $q\in \Lambda$ is
\begin{equation}\label{oca11}I(X|Y)(q):= \int \sum_{x=1}^d q(x) J(x*y) \log ( \frac{J(x*y) }{ \sum_{x^\prime} q(x^\prime) J(x^\prime*y)})d \mu(y).
\end{equation}

The capacity of $\mathcal{L}_J^*$ is 
$C=\sup \{ I(X|Y)(q) : q\in \Lambda\}$. In   Subsection \ref{anot} we analyze some examples where the non-zero capacity is explicitly computed.

\section{Non-Invariant Entropy (ni-entropy)} \label{sec:ni-entropy}

Let $\mathcal{X}_1 = \sigma^{-1} (\mathcal{F})$. For $j=1,2,\dots,d$ and $\rho\in \mathcal{M}$, the conditional expectation $E_{\rho} (I_{\overline{j}}| \mathcal{X}_1)$ is defined such that for any integrable function $b$:
\begin{equation} \label{luw} 
	\int (b \circ \sigma) E_{\rho} (I_{\overline{j}}| \mathcal{X}_1) \, d \rho = \int (b \circ \sigma) I_{\overline{j}} \, d \rho. 
\end{equation}

\begin{remark} \label{pola}  For each $j=1,2,...,d$, we choose a representative $E_{\rho} (I_{\overline{j}} \mid \mathcal{X}_1)$ which is well defined   and non-negative   in a set $K_j$  of $\rho$ positive measure, but not necessarily of $\rho$ measure  $1$ (see Example \ref{maili7}).  
$\sum_{j=1}^d E_{\rho} (I_{\overline{j}}| \mathcal{X}_1)=1$, $\rho$-a.e; not always $K_j\subset \overline{j}.$ For each $j$, $E_{\rho} ( I_{\overline{j}} | \mathcal{X}_1)=w_j \circ\sigma$ for  some integrable function $w_j$ (see Prop. 3.22 in \cite{Fer}).    Given $\rho \in \mathcal{M}$, for fixed $j$, $E_{\rho} ( I_{\overline{j}} | \mathcal{X}_1)>0$, $\rho$-a.e., $\Leftrightarrow$ there is no set $B\in \mathcal{X}_1$  such that  $I_{\overline{j}}=0$, $\rho$-a.e. in B (see Example \ref{moic16} where $E_{\rho} ( I_{\overline{3}} | \mathcal{X}_1) =0$, and $\log E_{\rho} ( I_{\overline{3}} | \mathcal{X}_1)$ is not integrable).
For $J \in \mathfrak{J}$,  the equilibrium probability $\mu_J$ satisfies  $E_{\mu_J} ( I_{\overline{j}} | \mathcal{X}_1) = J(j, \cdot)$ (see Example  \ref{maili}).
\end{remark}

\begin{definition} \label{ouent}
	The non-invariant entropy (ni-entropy) of $\rho \in \mathcal{M}$ is:
	\begin{equation} \label{feij4679} 
		h_{ni}(\rho) = - \int \sum_{j=1}^d E_{\rho} ( I_{\overline{j}} | \mathcal{X}_1) \log ( E_{\rho} ( I_{\overline{j}} | \mathcal{X}_1)) \, d \rho,
	\end{equation}
which is always finite.
	If for all $j=1,2,...,d$,  $\log ( E_{\rho} ( I_{\overline{j}} | \mathcal{X}_1))$ is $\rho$-integrable:
	\begin{equation} \label{filo} 
	h_{ni}(\rho) = - \int \sum_{j=1}^d I_{\overline{j}} \log ( E_{\rho} ( I_{\overline{j}} | \mathcal{X}_1)) \, d \rho.
	\end{equation}
\end{definition}
\begin{remark}\label{Ju}  $\forall j$,   $\log ( E_{\rho} ( I_{\overline{j}} | \mathcal{X}_1))$ is $\rho$-integrable $\Leftrightarrow E_{\rho} ( I_{\overline{j}} | \mathcal{X}_1)>0$, $\rho$ a. e. We refer the reader to Remark \ref{reret} to get a version of \eqref{filo} under no assumptions.
\end{remark}

It is known (see page 36 in \cite{PP}) that when $\mu$ is $\sigma$-invariant,  $\forall g\in \mathcal{C}$
\begin{equation} \label{lade} \sum_{j=1}^d\, \int g(j *x) \, E_{\mu} (I_{\overline{j}}| \mathcal{X}_1) (x)\, d \mu_1(x)=\int g\, d \mu.
\end{equation}

 \textbf{Crucially, this property may fail if $\mu$ is not $\sigma$-invariant} (see Example \ref{moj25}).   As for $J \in \mathfrak{J}$,  the equilibrium probability $\mu_J$ satisfies  $E_{\mu_J} ( I_{\overline{j}} | \mathcal{X}_1) = J(j, \cdot)$, we get that \eqref{lade} resembles the action of the dual operator $\mathcal{L}_J^*$.

Suppose a channel $\Phi$ transmits information. According to  page 49 in \cite{Feins1}:

{\it Shannon realized that while entropy is sufficient to characterize information content, it is not sufficient to characterize information transmission. Therefore, he defined a new quantity, conditional entropy, which measures the information lost to channel noise  during transmission of each element of a message.}

For $p \ll q,$ relative entropy $D(p||q)$ is also known as
Kullback-Leibler divergence: \,\,$D(p||q)= D((p_1,p_2,..,p_d)||(q_1,q_2,..,q_d))=  \sum_j p_j \log \frac{p_j}{q_j}  = $
\begin{equation} \label{opayu}\sum_{j=1}^d p_j \log p_j- \sum_{j=1}^d p_j \log q_j\geq 0.
\end{equation}

\noindent We refer to \cite{Vand} for this topic.
Here, a  natural concept of divergence is:

\begin{definition} \label{KL}
	The ni-Kullback-Leibler divergence $D_{ni}$ of $\rho$ given $\nu$ is: 
	$$0 \leq D_{ni}(\rho\,||\, \nu) =$$
$$ \int \sum_{j=1}^d E_{\rho} ( I_{\overline{j}} | \mathcal{X}_1)) \log ( E_{\rho} ( I_{\overline{j}} | \mathcal{X}_1)) \, d \rho - \int \sum_{j=1}^d E_{\rho} ( I_{\overline{j}} | \mathcal{X}_1)) \log ( E_{\nu} ( I_{\overline{j}} | \mathcal{X}_1)) \, d \rho,$$
in case the above expression is finite. As conditional expectations can vanish we can extend the values of  $D_{ni}(\rho\,||\, \nu)$ allowing  to be $+\infty$; for example, when $E_{\nu} ( I_{\overline{j}} | \mathcal{X}_1)) $ is zero on a set where $E_{\rho} ( I_{\overline{j}} | \mathcal{X}_1)) $ assigns positive mass. 
\end{definition}

\begin{remark} \label{Fra1}
Our definition of  $D_{ni}$ does not always coincide with the standard one, but it is  natural in our context.  Dissimilarities  occur  due to the fact that $E_{\rho} ( I_{\overline{j}} | \mathcal{X}_1))$ and $E_{\nu} ( I_{\overline{j}} | \mathcal{X}_1)$ may have different supports. 
For certain probabilities $\rho,\nu \in \mathcal{M}$, for example
when $\rho,\nu$ are on the image of $\mathcal{L}_J^*$, the property  $D_{ni}(\rho\,||\, \nu)=0$ does not  imply $\rho=\nu$ (see  Lemma \ref{ms1}). Expression (6.2.15) in \cite{Dembo} hints that $\log \int e^g d q$ has the status of a non-dynamical pressure $\mathfrak{p}(g)$.  Donsker-Varadhan variational  principle for  $D(p||q)$  claims:
\begin{equation} \label{plasx} \text{for}\,\, p\ll q, \,\,\,D(p||q)= \sup_{g}\{\int g\, dp - \log \int e^g d q\}=   \sup_{g}\{\int g\, dp - \mathfrak{p}(g)\}.
\end{equation}

In our setting, if $p,q\in \mathcal{M}$, it follows from \eqref{B33} that $D_{ni}(p,q)=$
\begin{equation} \label{kloe4} \sup_{g \in C^0} \{\int g dp - \mathcal{P}(g)  -\int\sum_j I_{\overline j}\log E_q(I_{\overline j}\mid\mathcal X_1)\,dp\}.\end{equation}
Comparing the two expressions above highlights how the  absolute continuity (always present in the classical case) should be appropriately handled here. $\mathcal{P}(g) $ does not mention $q$, and \eqref{kloe4} requires an extra term in $q$ (and $p$) to  match \ref{plasx}. The term $ \mathfrak{p}(\theta g)=\log \int e^{\theta g} d q$, $\theta \in \mathbb{R},$ in \eqref{plasx} - contains $q$ but not $p$ - it  is a form of cumulant generating function, and corresponds in \eqref{kloe4} to $ \mathcal{P}(\theta g)  +\int\sum_j I_{\overline j}\log E_q(I_{\overline j}\mid\mathcal X_1)\,dp$ (where no mention of $g$ in the last term). 
\end{remark}

 If $\rho,\nu\in \mathcal{G}$, and $D_{ni}(\rho\,||\, \nu)=0$, then $\rho=\nu$ (see Section \ref{inpro}).
The ni-Kullback-Leibler divergence 
$ D_{ni}(\mu_1\,||\, \mu_2) $  can be seen as a loss function in Decision Theory. In the case $\mu_1,\mu_2 \in \mathcal{G}$,  with respective Jacobians $J_1,J_2$, we get that $D_{ni}(\mu_1\,||\, \mu_2) =\int  \log J_1 \,d \mu_1 -\int \log J_2\, d \mu_1$, in line with Theorem $A$ in \cite{Ferre}, minimizing type II errors in Hypothesis Tests in a dynamical setting.

\begin{theorem} \label{iimp}
	The ni-entropy $h_{ni}(\rho)$ satisfies:
	\begin{enumerate}
		\item Upper semicontinuity and Borel measurability on $\mathcal{M}$.
		\item Concavity.
		\item Non-affinity.
		\item  $0 \leq h_{ni}(\rho) \leq \log d$, and the maximum can be attained outside $\mathcal{M}_\sigma$.
		\item Independence from the chosen version of $E_{\rho} ( I_{\overline{j}}| \mathcal{X}_1)$.
		\item Coincidence with the classical Kolmogorov-Sinai entropy for $\rho \in \mathcal{M}_\sigma$.
	\end{enumerate}
\end{theorem}

\begin{proof}
	
	First we  prove  1. and 2. It is known that for $t\geq 0$,  $t \log t = \sup_{u\in \mathbb{R}} ( u\, t - e^{u-1 })$. Denote by $\mathcal{W}$ the set of  $\mathcal{X}_1$-measurable continuous functions $\varphi:\Omega \to \mathbb{R}$, and for $j=1,2,...,d$, we set  $G_j := I_{\overline j} (x)\varphi -e^{ \varphi-1 }\in \mathcal{C}$; for a function $\varphi\in \mathcal{W}$, the value $\varphi(x)$ plays the role of the $u$ in the above $t \log t$.

	Given $j=1,2,...,d$, we claim that
	\begin{equation} \label{jiw}
	\int E_\rho(I_{\overline j}|\mathcal{X}_1) \log E_\rho(I_{\overline j}|\mathcal{X}_1)d \rho = 	\sup_{ \varphi \,\in \mathcal{W}} \{  \int I_{\overline j}(x) \varphi(x) -e^{ \varphi(x)-1 }\, d \rho(x)\}.
	\end{equation}
	
Note that for a given $j$
	\begin{equation} \label{jua}	E_\rho(I_{\overline j}|\mathcal{X}_1) \log E_\rho(I_{\overline j}|\mathcal{X}_1) =   	\sup_{ u \,\in \mathbb{R}}\{ E_\rho(I_{\overline j}|\mathcal{X}_1) u-e^{ u-1 }\}.
	\end{equation}
	
	For such  $j$, the values of all $\varphi \in \mathcal{W}$  can cover all  values $u\in \mathbb{R}$, thus
	\begin{equation} \label{juano}	E_\rho(I_{\overline j}|\mathcal{X}_1) \log E_\rho(I_{\overline j}|\mathcal{X}_1) =   	\sup_{ \varphi \,\in \mathcal{W}}\{ E_\rho(I_{\overline j}|\mathcal{X}_1)(x)  \varphi(x)-e^{ \varphi(x)-1 }\}.
	\end{equation}

Using \eqref{juano}, and the property that for a $\mathcal{X}_1$-measurable function  $\varphi$
$$\int E_\rho(I_{\overline j}|\mathcal{X}_1) (x)\varphi(x) d \rho(x)= \int I_{\overline j}(x) \varphi(x) d \rho(x), \,\,\,\text{then, we get}$$
$$	\int E_\rho(I_{\overline j}|\mathcal{X}_1) \log E_\rho(I_{\overline j}|\mathcal{X}_1)d \rho  \stackrel{1}{=}  	\sup_{ \varphi \,\in \mathcal{W}}\{\int E_\rho(I_{\overline j}|\mathcal{X}_1)(x) \varphi(x) -e^{ \varphi(x)-1 }d \rho(x)\}=$$
$$ 	\sup_{  \varphi \,\in \mathcal{W}}\{\int I_{\overline j}(x) \varphi(x) -e^{ \varphi(x)-1 }d \rho(x)\}= \sup_{  \varphi \,\in \mathcal{W}}\{ \int G_j(x)d\rho(x)\},$$
showing \eqref{jiw}.
 The interchange of  supremum and integral used above in  the transition within the equality $ \stackrel{1}{=}$ requires the use of the following lemma:

\begin{lemma} \label{req}
	For a given $j=1,2,...,d$,
	\begin{equation} \label{jua15}\sup_{ \varphi \,\in \mathcal{W}}	\{	\int  E_\rho(I_{\overline j}|\mathcal{X}_1) \varphi -e^{ \varphi-1 }d \rho\}= 	\int \sup_{ \varphi \,\in \mathcal{W}}\{ E_\rho(I_{\overline j}|\mathcal{X}_1) \varphi -e^{ \varphi-1 }\} d \rho.
	\end{equation}
	
	\end{lemma}
\begin{proof} We will use the theory of normal integrands (Theorem 16.40  in \cite{Rock}).
	
	1). For each $u$, we get that $x \to E_\rho(I_{\overline j}|\mathcal{X}_1) (x)u-e^{ u-1 }$ is $\rho$-integrable.
	
	2) for each $x$ the function $ u \to  E_\rho(I_{\overline j}|\mathcal{X}_1)(x) u-e^{ u-1 }$ is continuous.
	
	This shows that $(x,u)\to E_\rho(I_{\overline j}|\mathcal{X}_1) (x)u-e^{ u-1 }$ is  normal. The space $\mathcal{W}$ is a decomposable Banach space and  $x \to E_\rho(I_{\overline j}|\mathcal{X}_1) (x)\varphi(x)-e^{ \varphi(x)-1 }$ is always $\rho$-integrable.  \eqref{jua15} follows from Theorem 3.4   \cite{Garr}, or Proposition 2.1  \cite{Mor}.\end{proof}

As,  $\rho \to \int E_\rho(I_{\overline j}|\mathcal{X}_1)  \log E_\rho(I_{\overline j}|\mathcal{X}_1)d \rho$ is the supremum of a collection of linear continuous operators, then it is  lower semicontinuous. If $\nu_k\to \rho$, then
\begin{equation}\label{ufer}\int E_\rho(I_{\overline j}|\mathcal{X}_1) \log E_\rho(I_{\overline j}|\mathcal{X}_1)d \rho \leq \liminf_{k \to \infty} 
\int E_{\nu_k}(I_{\overline j}|\mathcal{X}_1) \log E_{\nu_k}(I_{\overline j}|\mathcal{X}_1)d \nu_k.
\end{equation}

 Therefore,  $-h_{ni}(\rho) =\int \sum_{j=1}^d E_\rho(I_{\overline j}|\mathcal{X}_1) \log E_\rho(I_{\overline j}|\mathcal{X}_1)d \rho \leq$
 \begin{equation}\label{bufer} \liminf_{k \to \infty} 
 \int \sum_{j=1}^d E_{\nu_k}(I_{\overline j}|\mathcal{X}_1) \log E_{\nu_k}(I_{\overline j}|\mathcal{X}_1)d \nu_k= \lim_{k \to \infty} - h_{ni}(\nu_k).
 \end{equation}
	
	 It follows that $\rho \to  h_{ni}(\rho)$ is upper semicontinuous and concave.

	\medskip

	3. Consider in Example \ref{mojic}, for $\lambda\in [0,1]$, the family of probabilities $\mu^\lambda$, where $p_1=\lambda$ and $p_2=(1-\lambda).$
	The function $\mu^\lambda \to h_{ni}(\mu^\lambda)$ is not affine.
	
	\medskip
	
	4. As $\sum_{j=1}^d I_{\overline{j}}=1,$ we get that $\sum_{j=1}^d E_{\rho} (I_{\overline{j}}| \mathcal{X}_1)= 1,$ $\rho$-a.e. Thus,  $E_{\rho} (I_{\overline{j}}| \mathcal{X}_1)\leq 1.$ Given $j$,
	$   E_{\rho} (I_{\overline{j}}| \mathcal{X}_1)$ is nonnegative $\rho$-a.e because $I_{\overline{j}}\geq 0.$ Moreover, the function $x \log x$ is nonpositive for $x\in [0,1].$ Then, $h_{ni}(\rho)\geq 0.$ As $- \sum_{j=1}^d p_j \log p_j \leq \log d$,  $\forall p=(p_1,p_2,...,p_d)$, we get for $\rho$-a.e $x\in \Omega$, $-\sum_{j=1}^d  E_{\rho}( I_{\overline{j}} | \mathcal{X}_1)(x)  \log ( E_{\rho} ( I_{\overline{j}} | \mathcal{X}_1))(x)\leq  \log d$, and so $h_{ni}(\rho)\leq \log d.$
	
	In Example \ref{mojic} we exhibit  $\mu\notin \mathcal{M}_\sigma $ satisfying $h_{ni}(\mu)=\log d.$

	\medskip
	
	5.  The value of the entropy $h_{ni}(\rho)$ does not depend on the choice of  $E_{\rho} ( I_{\overline{j}}| \mathcal{X}_1)$, $j=1,2,...,d$, for a general $\rho \in \mathcal{M}$. This is so  because for fixed $j$ we have to integrate over the support $\rho_0$, and $E_{\rho} ( I_{\overline{j}}| \mathcal{X}_1)$ is unique $\rho$-a.e.

	\medskip

	To prove 6., we refer the reader to Section 3 in \cite{PP} (see also Section 8 in \cite{PY}):   when $\mu  \in \mathcal{M}_\sigma$,  \eqref{feij4679}   is the  Kolmogorov-Sinai entropy.
	\end{proof}

The  operator $L_1: \mathcal{C}^+\to \mathcal{C}^+$  satisfies $L_1 (u)(x) =v(x)=  \sum_{j=1}^d u(j,x)$. 

\begin{theorem} \label{tirri}  Given $\rho \in \mathcal{M}$, then 
	\begin{equation} \label{este59}
	h_{ni}(\rho) = \inf_{u \in \mathcal{C}^+} \int \log \left( \frac{L_1 (u) \circ \sigma}{u} \right) d\rho.
	\end{equation}
	
\end{theorem}

We need to prove a Lemma before going to  the proof of Theorem \ref{tirri}.

\begin{lemma} \label{ecom} Given $\rho \in \mathcal{M}$, there exists a sequence of functions $c_n^j \in \mathcal{C}^+$ which are $\mathcal{X}_1$-measurable,  $j=1,2,...,d$, $n\ \in \mathbb{N}$, satisfying 
\begin{equation} \label{qwr} \lim_{n \to \infty} \int \sum_{j=1}^d c_n^j  \log c_n^j  d \rho=\int \sum_{j=1}^d E_{\rho} (I_{\overline{j}} \mid \mathcal{X}_1) \log (E_{\rho} (I_{\overline{j}} \mid \mathcal{X}_1)\,)  d\rho.\end{equation}
	
\end{lemma}
\begin{proof}
	Given $j=1,2,...,d$,  the $\mathcal{X}_1$-measurable function $\,E_{\rho} (I_{\overline{j}}  \mid \mathcal{X}_1):\Omega \to [0,1]$ may achieve the value $1$ or $0$, in the support of $\rho$. Via Lusin's theorem and the Tietze extension theorem,  given $k,n$, one can show that  $ \,E_{\rho} (I_{\overline{j}}  \mid \mathcal{X}_1) $ can be approximated by $\frac{1}{2 d n}$ in $L^1 (\rho)$ by a continuous $\mathcal{X}_1$-measurable function $c^{j}_n:\Omega \to (0,1]$, s. t., denoting $K_{j,n}=\{x : c^{j}_n (x)-  E_{\rho}(I_{\overline{j}}  \mid \mathcal{X}_1)(x)<0\}$, it is valid $\rho(K_{j,n})<\frac{1}{4 d n}$. Note that $ |\, c_n^{j} -E_{\rho} (I_{\overline{j}} \mid \mathcal{X}_1) \,|\leq 1.$ $\int |c_n^{j} -E_{\rho} (I_{\overline{j}} \mid \mathcal{X}_1) \,|\, d \rho=$
	\begin{equation} \label{flule}\int I_{\Omega -K_{j,n}}( c_n^{j} -E_{\rho} (I_{\overline{j}} \mid \mathcal{X}_1) \,)\, d \rho + \int I_{K_{j,n}}|\, c_n^{j} -E_{\rho} (I_{\overline{j}} \mid \mathcal{X}_1) \,|\, d \rho  \leq \frac{1}{d n}.\end{equation}
	
	It is known that $\sup_{x\in[0,1]} |y \log y |=  e^{-1}.$ Then, for all $x\in \Omega$
	\begin{equation} \label{flul}|c_n^j  (x)\log c_n^j(x)  -   E_{\rho} (I_{\overline{j}} \mid \mathcal{X}_1)(x) \log (E_{\rho}(I_{\overline{j}} \mid \mathcal{X}_1))(x)\, |\leq 2 e^{-1}<2.\end{equation}
	For $0\leq s\leq  r\leq 1$, it is valid $r-s \geq r\log r - s \log s$. Given $j$, taking $r=c_n^j (x) \geq 0$ and $s= E_{\rho} (I_{\overline{j}} \mid \mathcal{X}_1) (x)$, $x \in \Omega - K_{j,n}$, and integrating on $\rho$  
	$$\int| c_n^j  \log c_n^j   -   E_{\rho} (I_{\overline{j}} \mid \mathcal{X}_1) \log (E_{\rho}(I_{\overline{j}} \mid \mathcal{X}_1))|d \rho \,  \leq$$
	$$\int I_{  K_{j,n} }\,|\, c_n^j  \log c_n^j   -  E_{\rho} (I_{\overline{j}} \mid \mathcal{X}_1) \log (E_{\rho}(I_{\overline{j}} \mid \mathcal{X}_1)\,|
	d \rho \,\,+$$
	$$\int I_{  \Omega - K_{j,n} }\,|c_n^j  \log c_n^j   -  E_{\rho} (I_{\overline{j}} \mid \mathcal{X}_1) \log (E_{\rho}(I_{\overline{j}} \mid \mathcal{X}_1)|
	d \rho \,\,\leq$$
	$$ \frac{1}{2 d n} + \int  I_{  \Omega - K_{j,n} }|\,c_n^j  -   E_{\rho} (I_{\overline{j}} \mid \mathcal{X}_1)\,|d \rho\leq$$
	\begin{equation} \label{klead}\frac{1}{2 d n} + \int  |\,c_n^j  -   E_{\rho} (I_{\overline{j}} \mid \mathcal{X}_1)\,|d \rho\leq \frac{1}{d n} .
	\end{equation}
	Finally,  \eqref{qwr} follows from the existence of the family $c_n^j $, $j=1,2,...,d$, s. t.
	\begin{equation} \label{klead8}\int |\sum_j c_n^j  \log (c_n^j)   - \sum_j   E_{\rho} (I_{\overline{j}} \mid \mathcal{X}_1|) \log (E_{\rho}(I_{\overline{j}} \mid \mathcal{X}_1))\,|d \rho\leq  \frac{1}{  n}.
	\end{equation}
	
	Given  $a=1,2,...,d$, due to \eqref{flule}
	$$ \int |\,(c^{j}_n (a,x_2,...)- E_{\rho} (I_{\overline{j}}  \mid \mathcal{X}_1)(a,x_2,...)\,\,| d \rho(x)\leq $$
	\begin{equation} \label{fluiu} \int  |c_n^{j} -E_{\rho} (I_{\overline{j}} \mid \mathcal{X}_1) \,|\, d \rho  \leq \frac{1}{d n}.
	\end{equation}
	Moreover, given  $a=1,2,...,d$,   from \eqref{klead}
	$$\int| c_n^j (a,\sigma(x))  \log c_n^j   (a,\sigma(x))-  $$
	$$ E_{\rho} (I_{\overline{j}} \mid \mathcal{X}_1)  (a,\sigma(x))\log (E_{\rho}(I_{\overline{j}} \mid \mathcal{X}_1)) (a,\sigma(x))|d \rho(x) \,  \leq$$
	\begin{equation} \label{fluiu3}\int | c_n^j (x)  \log c_n^j   (x)-   E_{\rho} (I_{\overline{j}} \mid \mathcal{X}_1)  (x)\log (E_{\rho}(I_{\overline{j}} \mid \mathcal{X}_1)) (x)|d \rho(x) \,  \leq\frac{1}{d n}.
	\end{equation}
\end{proof}
\medskip

\begin{remark} \label{reret} Given $\rho$, $\forall j$, we set  $U_j= \{ E_{\rho} (I_{\overline{j}} \mid \mathcal{X}_1) >0\}$, which is a set  in $\mathcal{X}_1$. Denote $V_j(x) =I_{U_j} (x)E_{\rho} (I_{\overline{j}} \mid \mathcal{X}_1)(x)$,  which is $\mathcal{X}_1$-measurable. Then,   1) $V_j$ is $\rho$-integrable. 2) $\int V_j \log (V_j) d \rho = \int E_{\rho} (I_{\overline{j}} \mid \mathcal{X}_1) \log (E_{\rho} (I_{\overline{j}} \mid \mathcal{X}_1)) d \rho.$ 3) $ \int I_{\overline{j}} \log (E_{\rho} (I_{\overline{j}} \mid \mathcal{X}_1)) d \rho= \int E_{\rho} (I_{\overline{j}} \mid \mathcal{X}_1) \log (E_{\rho} (I_{\overline{j}} \mid \mathcal{X}_1)) d \rho.$ Therefore, $h_{ni}(\rho) =- \sum_{j=1}^d\int V_j \log (V_j) d \rho$.
	
	\end{remark}

Proof of Theorem \ref{tirri}.

\begin{proof} 
	
	For any $u \in \mathcal{C}^+$, let $q_j(x) = \frac{u(j ,\sigma(x) )}{(\mathcal{L}_1u)\,(\sigma(x))}>0$. By definition, $\sum_j q_j(x) = 1$.   Denote $p_j (x)= E_{\rho}(I_{\overline{j}} | \mathcal{X}_1)(x)$, now using the classical inequality  
	\begin{equation} \label{opayu}\sum_{j=1}^d p_j \log p_j- \sum_{j=1}^d p_j \log q_j\geq 0,
	\end{equation}
	integrating both sides with respect to $\rho$, from the expected value property:
	$$  h_{ni}(\rho) \leq - \int \sum_{j=1}^d E_{\rho}(I_{\overline{j}} | \mathcal{X}_1)(x) \log \left( \frac{u(j, \sigma(x))}{L_1 (u)( \sigma(x))} \right) d\rho (x) = $$
	\begin{equation} \label{hum1}   - \int \sum_{j=1}^d I_{\overline{j}}  \log \left( \frac{u(j, \sigma(x))}{L_1(u)(\sigma(x))} \right) d\rho(x)=  \int \log \left( \frac{L_1 (u)(\sigma( x))}{u(x)} \right) d\rho(x).
	\end{equation}

	If we take the infimum over all $u \in \mathcal{C}^+$ we get
	\begin{equation} \label{lun2}h_{ni}(\rho) \leq  \inf_{u \in \mathcal{C}^+} \int \log \left( \frac{L_1 (u) \circ \sigma}{u} \right) d\rho.
	\end{equation}
	
	For the other inequality $\geq$ in \eqref{este5}, we will use the notation and the claims of Remark \ref{reret}.   Take 
	 $v=e^{\sum_{j=1}^d V_j \log (V_j) }$. Note that $\sum_{j=1}^d V_j \log (V_j)$ is not continuous, but we follow our reasoning with this caveat ($v$ will be approximated by  $u_n\in \mathcal{C}^+$,  and the infimum in \eqref{este59} will not be attained). It follows
	 $$\int \log (v)\,  d \rho =\int \sum_{j=1}^d  E_{\rho} (I_{\overline{j}} | \mathcal{X}_1) \log ( E_{\rho} ( I_{\overline{j}} | \mathcal{X}_1))d \rho = - h_{ni}(\rho).$$
\eqref{klead} shows the existence of a sequence of continuous functions $c_n^j$, $n \in \mathbb{N}$, $j=1,2,...,d$, such that $\int \sum_j c_n^j\log c_n^j d \rho$ approximates the value $- h_{ni}(\rho).$ Denote $u_n (x)= e^{ \sum_j c_n^j(x_1,x_2,...)\log c_n^j(x_1,x_2,...)}$,  $n \in \mathbb{N}$. Then, $\lim_{n \to \infty}\int \log u_n d \rho =- h_{ni}(\rho)$.  We will  show that $\int \log \left( L_1 (u_n)\right) d\rho\leq 0$ for  some large $n$. Note that
	\begin{equation} \label{luna23} L_1 (u_n) (\sigma(x)) =\sum_{a=1}^d e^{ \sum_{j=1}^d c_n^j(a,\sigma(x))\log c_n^j (a,\sigma(x))}.\end{equation}
	and
	\begin{equation} \label{fluiu5}L_1 (v) (\sigma(x)) =\sum_{a=1}^d e^{\sum_{j=1}^d V_j (a,\sigma(x))\log (V_j)(a, \sigma(x))  }.
	\end{equation}

	We want to show first  that 
	$ \int\log(\,L_1 (v)(\sigma(x))) d \rho (x)\leq 0.$

	By the property that a $\mathcal{X}_1$-measurable function does not depend on the first variable $x_1$, Jensen inequality  and \eqref{fluiu5}, we get

$$ \int\log(\,L_1 (v)(\sigma(x))) d \rho (x)\leq \log( \int L_1 (v)(\sigma(x))) d \rho (x))= $$
$$ \log (\int \sum_{a=1}^d e^{\sum_{j=1}^d V_j (a,\sigma(x))\log (V_j)(a, \sigma(x))  }     d \rho (x))=$$
$$ \log( \int  \sum_{a=1}^d e^{\sum_{j=1}^d I_{\overline{j}} (a,\sigma(x))\log (I_{U_j}(a, \sigma(x))  E_{\rho} ( I_{\overline{j}} | \mathcal{X}_1)) (a, \sigma(x)))}     d \rho (x))=$$
$$ \log( \int  \sum_{a=1}^d e^{ I_{\overline{a}} (a,\sigma(x))\log  (I_{U_a}(a, \sigma(x)) E_{\rho} ( I_{\overline{a}} | \mathcal{X}_1))  (a, \sigma(x))}     d \rho (x))=$$
$$ \log( \int  \sum_{a=1}^d e^{\log  (I_{U_a}(a, \sigma(x))  E_{\rho} ( I_{\overline{a}} | \mathcal{X}_1))  }    d \rho (x))\stackrel{2}{=}$$
$$\log( \int  \sum_{a=1}^d I_{U_a}(a, \sigma(x))  E_{\rho} ( I_{\overline{a}} | \mathcal{X}_1)  (a, \sigma(x))    d \rho (x))=$$
$$\log( \int  \sum_{a=1}^d I_{U_a}(x)  E_{\rho} ( I_{\overline{a}} | \mathcal{X}_1)(x)    d \rho (x))=$$
	\begin{equation} \label{fluiu6} \log( \int  \sum_{a=1}^d  E_{\rho} ( I_{\overline{a}} | \mathcal{X}_1)  (x)    d \rho (x))=0.\end{equation}
	
Passage $\stackrel{2}{=}$ requires to split the integration on $\rho$ when $ x\in U_a $ or $x\notin U_a$.
	
	Then, $  \log(L_1 (v)(\sigma(.)) $ is integrable. Moreover, given $a$, from \eqref{fluiu3},  in $L^1(\rho)$
	$$ -e^{-1} d \leq\sum_{j=1}^d c_n^j(a,\sigma(.))\log c_n^j (a,\sigma(.)) \to $$
	\begin{equation} \label{fluiu4} \sum_{j=1}^d     E_{\rho} ( I_{\overline{j}} | \mathcal{X}_1))(a,\sigma(.)\log ( E_{\rho} ( I_{\overline{j}} | \mathcal{X}_1))(a,\sigma(.))\leq 0.\end{equation}
	
\begin{remark}\label{doik}	When $ -e^{-1} d\leq  y_a \leq 0$, for all $a=1,2...,d$, 
	the function $G(y_1,y_2,...,y_d)=\log (\sum_{a=1}^d e^{y_a})$ is Lipschitz. For each $a=1,2...,d$, $n \in \mathbb{N},$ consider  a sequence   of continuous functions $ f_a^n : \Omega \to  [- e^{-1}d,0]$. Assume, that  in $L^1(\rho)$,  the sequence  $f_a^n  \to f_a$, when $n \to \infty$, then, it follows that in $L^1(\rho)$, $\log(\sum_{a=1}^d e^{f_a^n})  \to \log(\sum_{a=1}^d e^{f_a} )$, when $n \to \infty$.
	\end{remark}

 Set $f_a(.)= \sum_{j=1}^d E_{\rho} ( I_{\overline{j}} | \mathcal{X}_1))(a,\sigma(.)\log ( E_{\rho} ( I_{\overline{j}} | \mathcal{X}_1))(a,\sigma(.))$,  and  $f_a^n(.)=  \sum_{j=1}^d c_n^j(a,\sigma(.))\log c_n^j (a,\sigma(.))$, from Remark \ref{doik}, \eqref{luna23},\eqref{fluiu5}, \eqref{fluiu6} and\eqref {fluiu4} it follows that given $\epsilon>0$, we get that $\int \log \left( L_1 (u_n)\right) d\rho\leq \epsilon,$  for large $n$. Then,
	$$h_{ni}(\rho) \geq  \inf_{u \in \mathcal{C}^+} \int \log \left( \frac{L_1 (u)\circ \sigma}{u} \right) d\rho.$$
	From the above,  \eqref{este59} follows.\end{proof}

\medskip

\begin{remark} \label{klio}
	Suppose $\mu_3= \mathcal{L}^*_{J} (\mu_1)$. For any sufficiently small open set $A$ around $y$, $\mu_3(A)>0$ if and only if $\mu_1 (\sigma^{-1} (A))>0$.
From $\mu_3=\mathcal L_J^*(\mu_1)$, for any cylinder $C$, one has $\mu_3(C)=\int\sum_iJ(i,x)I_C(i*x)\,d\mu_1(x)$ (see Example \ref{moj0}).   Due to the strict positivity of $J$, the above equivalent positivity relation  for $C$ is governed by the tails $x$ having at least one branch $i*x\in C$. The cylinder sets generate the Borel sigma-algebra. Given the property that a probability can be described by inner and outer measures \cite{Fer}, it follows that the equivalence holds for open sets $A$.  Therefore $\operatorname{supp}\mu_3=\sigma^{-1}(\operatorname{supp}\mu_1)$.

\end{remark}

\begin{proposition} \label{cont}
	Let $J$ be the H\"older Jacobian of  $\mu\in \mathcal{G}$. For $\mu_1\in \mathcal{M}$, let $\mu_3= \mathcal{L}^*_{J} (\mu_1)$.
Then, the conditional expectation is:
	\begin{equation} \label{moj8}
		E_{\mu_3} ( I_{\overline{j}} | \mathcal{X}_1) (x_1,x_2,\dots)= J (j, x_2,x_3,\dots)>0,\,\, \mu_3 \,\, \text{almost everywhere}.
	\end{equation}
	
	Moreover, the entropy formula linking $\mu_1$ and $\mu_3$ is:
	\begin{equation} \label{moj47}
		h_{ni}(\mu_3) = - \sum_{j=1}^d \int J (j, x) \log J (j, x)\, d \mu_1 (x).
	\end{equation}
	Alternatively, using $\mathfrak{T} (\mu_3)=\mu_1$, we have:
	\begin{equation} \label{moj47x }
		h_{ni}(\mu_3) = - \sum_{j=1}^d \int J (j, \sigma(x)) \log J (j, \sigma(x)) d \mu_3 (x).
	\end{equation}
	Finally, $\mu_3$ is $\sigma$-invariant if and only if $\mu_1=\mu$ (in which case $\mu_3=\mu$).
\end{proposition}

\begin{proof}
	For a fixed $j$, there exists a measurable $w_j$ such that $E_{\mu_3} ( I_{\overline{j}} | \mathcal{X}_1)= w_j \circ \sigma$. For any measurable function $g$:
	$$\int I_{\overline{j}} (g \circ \sigma) d \mu_3 = \int \mathcal{L}_{J} ( I_{\overline{j}} (g \circ \sigma)) d \mu_1 
		= \int J (j, x) g (x) d \mu_1 (x).$$
Simultaneously, by definition of conditional expectation:
	$$ \int E_{\mu_3} ( I_{\overline{j}} | \mathcal{X}_1) (g \circ \sigma) d \mu_3 = \int (w_j \circ \sigma)(g \circ \sigma) d \mu_3 = \int w_j g d \mu_1.$$
	Equating both sides, we can choose $w_j (x)= J (j, x)$  in the support of $\mu_1$. From Remark \ref{klio}, given $f_1,f_2:\Omega \to \mathbb{R}$, 
if  $f_1=f_2$ in the support of $\mu_1$, then  $f_1 \circ \sigma= f_2 \circ  \sigma$ in the  support of $\mu_3$, establishing \eqref{moj8}.
	
	Using \eqref{moj8} and the duality $\int f d\mu_3 = \int \mathcal{L}_J(f) d\mu_1$ applied to the function  $f(x) = I_{\overline{j}}(x) \log J(j, \sigma(x))$, we prove \eqref{moj47}:
	$$h_{ni}(\mu_3) = - \sum_{j=1}^d \int \mathcal{L}_{J}\big( I_{\overline{j}} \cdot (\log J(j, \cdot) \circ \sigma) \big) (x) d \mu_1 (x)= $$
		$$ - \sum_{j=1}^d \int \log J (j, x) \mathcal{L}_{J}(I_{\overline{j}})(x) d \mu_1 (x) 
		= - \sum_{j=1}^d \int J(j, x) \log J (j, x) d \mu_1 (x).$$
	
	Assume $\mu_3$ is $\sigma$-invariant. Using \eqref{lade} and \eqref{moj8}, for any continuous $g$:
	$$ \int g d \mu_3 = \sum_{j=1}^d \int g(j,\sigma(x)) J(j,\sigma(x)) d \mu_3(x) = \int \mathcal{L}_{J} (g) d \mu_3 = \int g d\mathcal{L}_{J}^* (\mu_3).$$
	Thus $\mathcal{L}_{J}^* (\mu_3)=\mu_3$, which implies $\mu_3=\mu$ (since $\mu$ is the unique fixed probability). Consequently, $\int g d \mu_1 = \int (g \circ \sigma) d \mathcal{L}_{J}^*(\mu_1) = \int (g \circ \sigma) d \mu_3 = \int (g \circ \sigma) d \mu = \int g d \mu$, ensuring $\mu_1=\mu$.
\end{proof}

\begin{proposition} \label{moj3}
	Given $\mu_1 \in \mathcal{M}$ and $\mu_3= \mathcal{L}^*_{J} (\mu_1)$, it follows that:
	\begin{equation} \label{jui}
		h_{ni}(\mu_3)= - \int \log J d \mu_3 > 0.
	\end{equation}
\end{proposition}

\begin{proof}
	Applying \eqref{moj8} and the transfer operator properties: $h_{ni}(\mu_3) =$
	$$	 -\sum_{j=1}^d \int I_{\overline{j}} \log ( E_{\mu_3} ( I_{\overline{j}} | \mathcal{X}_1)) d \mu_3 = -\sum_{j=1}^d \int \mathcal{L}_{J} \big(I_{\overline{j}} \cdot (\log J(j,\cdot) \circ \sigma)\big) d \mu_1 =$$
		$$ -\sum_{j=1}^d \int J(j, x) \log J(j,x) d \mu_1 (x) =
		 -\int \mathcal{L}_{J} (\log J) d \mu_1 = - \int \log J d \mu_3.$$\end{proof}

In the next example the Jacobian is recovered from the expected value formalism when considering H\"older Gibbs probabilities.

\begin{example} \label{maili}
	Suppose for all $j$,  $x \to J(j\,x)>0$ are H\"older continuous functions such that $\sum_{j=1}^d J(j\,x)=1$, for all $x\in \Omega$. $\mu_J$ denotes the Gibbs probability for $J\in \mathfrak{J}$. It is known (see \cite{PP}) that for all $g\in \mathcal{C}$, $E_{\mu_J} (g| \mathcal{X}_1) = \mathcal{L}_{J} (g) \circ \sigma.$
	In this way, for any $j$ and $x=(x_1,x_2,x_3,\dots)$, we have
	$$E_{\mu_J} (I_{\overline{j}}| \mathcal{X}_1)(x_1,x_2,x_3,\dots) = \mathcal{L}_{J} (I_{\overline{j}} ) \circ \sigma(x_1,x_2,x_3,\dots) =
	J(j,x_2,x_3,\dots)= $$
	$$ J(j,\sigma(x_1,x_2,x_3,\dots)) 
	= (w_j \circ \sigma)(x_1,x_2,x_3,\dots),$$
	where $w_j(x_1,x_2,x_3,\dots) = J(j,x_1,x_2,x_3,\dots)$ (see Prop. 3.4 in \cite{PP}). Then, 
	$$ h(\mu_J)= h_{ni}(\mu_J) =  - \int \sum_{j=1}^d E_{\mu_J} ( I_{\overline{j}} | \mathcal{X}_1)\log ( E_{\mu_J} ( I_{\overline{j}} | \mathcal{X}_1)) d \mu_J =$$
	$$ - \int \sum_{j=1}^d I_{\overline{j}} \log ( E_{\mu_J} ( I_{\overline{j}} | \mathcal{X}_1)) d \mu_J = -  \int \sum_{j=1}^d I_{\overline{j}}(j,\sigma(x)) \log J(j,\sigma(x)) d\mu_J(x)=$$
	\begin{equation} \label{lort}  -\sum_{j=1}^d \int  \log J(j,\sigma(x)) d\mu_J(x)=- \int \log J d \mu_J.
	\end{equation}
	
\end{example}
\section{Information projections} \label{inpro}

In this Section, our focus is on Information Projections associated with the ni-Kullback-Leibler divergence.
Consider, for instance,  the variational  problem:   fix a probability $\mu_2\in \mathcal{G}$,  $f:\Omega \to \mathbb{R}$ and $y\in \mathbb{R}$, and find the minimizing $\rho$ for
\begin{equation} \label{tye} D_y=\inf_{\rho \in \mathcal{M}}\left\{D_{ni} (\rho,\mu_2) \,|\, \int f d \rho=y  \right\}\geq 0.\end{equation}

The constraint set can be empty.  Set  $A_f:=\{y | \Pi_{y} \neq \emptyset\}=[\min f, \max f]$.    From now on we will just consider values $y\in A_f.$   
In the classical setting (non-dynamical and also not within our setting) the probability $\rho$ attaining the minimal value in \eqref{tye} is called an information projection (see also \cite{LR3}); it can sometimes be expressed as a tilting probability (see  Lemma 3 in \cite{LM34}, or (3.22) in \cite{Csi1}); here we get similar (but not exactly the same) kinds of results.

Assume the H\"older Jacobians $J_1,J_2 \in\mathfrak{J}$ correspond, respectively, to the equilibrium probabilities $\mu_1,\mu_2\in \mathcal{G}$. As  $\mu_1 \in \mathcal{M}_\sigma$, we get from Example \ref{maili}
$$ \int E_{\mu_1} (I_{\overline{j}}| \mathcal{X}_1)\log ( E_{\mu_1} (I_{\overline{j}}| \mathcal{X}_1))- \int E_{\mu_1} (I_{\overline{j}}| \mathcal{X}_1)\log ( E_{\mu_2} (I_{\overline{j}}| \mathcal{X}_1))=$$
$$ \int J_1(j,\sigma(x)) \log J_1(j,\sigma(x)) d \mu_1 (x)-\int J_1(j,\sigma(x)) \log J_2(j,\sigma(x)) d \mu_1(x)=   $$
$$ \int J_1(j,x) \log J_1(j,x) d \mu_1 (x)-\int J_1(j,x) \log J_2(j,x) d \mu_1(x).   $$

Now, summing over $j$ the last expression, as $\mathcal{L}_{J_1}^* (\mu_1)=\mu_1$, we get
$$D_{ni} (\mu_1 \| \mu_2)= \int  \log J_1 d \mu_1 -\int \log J_2 d \mu_1 =-h(\mu_1) -\int \log J_2 d \mu_1 .   
$$

One can show that if $\mu_1,\mu_2 \in \mathcal{G},$ then $D_{ni} (\mu_1 \| \mu_2)=0$ implies $\mu_1= \mu_2$. This follows from the uniqueness of the equilibrium state for $P(\log J_2)=0.$

We will be interested in the following particular problem:
given $y\in \mathbb{R}$, and a H\"older function $f:\Omega \to \mathbb{R}$,  consider a variable $\mu_1\in \mathcal{M}_\sigma$, and then set
$$\Pi_{y}^\sigma= \{\mu_1 \in \mathcal{M}_\sigma\,|\, \int f d \mu_1=y \},\,\,\text{which is a convex set}.$$ 

We fix $\mu_2\in \mathcal{G}$.
Denote $D_y^\sigma= D_{ni} (\Pi_y^\sigma, \mu_2)= \inf \{D_{ni}(\mu_1, \mu_2)\,|\,\mu_1\in \Pi_y^\sigma\}.$ 
A probability  $\mu_1$ minimizing $ D_y^\sigma $ will be called the dynamical ni-Information Projection of the triple $\mu_2$, $f$ and $y$ (see  \cite{LR3} for related results).

An equivalent problem to minimizing $ D_y^\sigma $ is: find $\mu_1$ maximizing
$$\sup \{ h(\mu_1) +\int \log J_2 d \mu_1:\\ \mu_1 \in \Pi_{y}^\sigma\}. $$

\begin{theorem} Given $y$, assume that $\lambda_y$ satisfies $\frac{d}{d \lambda} P(\log J_2 + \lambda f)|_{\lambda_y}=y$.
	Denote by $\mu^{\lambda_y}$ the equilibrium probability for $P( \log J_2 + \lambda_y f )$, and by $h_{\lambda_y}$ the eigenfunction for the Ruelle operator $\mathcal{L}_{\log J_2 + \lambda_y f}$. Then, $\mu^{\lambda_y}$ is the  dynamical ni-Information projection for  $\mu_2$, $f$, $y$; and  for any continuous function $\psi$:
	\begin{equation} \label{triuvi} \lim_{k \to \infty} \int \, \psi(x)\,\frac{e^{\,\sum_{j=0}^{k-1} \lambda_y f (\sigma^j (x))}}{ \int e^{\, \sum_{j=0}^{k-1}  \lambda_y f  (\sigma^j (z) )}d \mu_{\log J_2}(z)}h_{\lambda_y}(x) d \mu_{\log J_2} (x) =\int \psi  d \mu^{\lambda_y}.
	\end{equation} 
\end{theorem}

\begin{proof}
	
	Consider the following problem: maximize over $\mu_1\in \mathcal{M}_\sigma$
	\begin{equation} \label{reri}G(\mu_1) = h(\mu_1) +\int \log J_2 d \mu_1+ \lambda \int f d \mu_1= h(\mu_1) +\int (\log J_2 + \lambda  f) d \mu_1. \end{equation}
 In this direction,  we have to analyze  the pressure problem $P(\log J_2 + \lambda f)$. Denote by $\mu^\lambda $ the corresponding equilibrium probability. As we assumed that $\log J_2$ and $f$ are H\"older, we get that  $\mu^\lambda \in\mathcal{G}$.
	Assume the existence of $\lambda_y$ such that $\frac{d}{d \lambda} P(\log J_2 + \lambda f)|_{\lambda_y}=y$.  The values that can be reached by such a derivative are in an open finite interval (this follows from \cite{PP}). From \cite{PP} we also get that $y= \int f d \mu^{\lambda_y}$.   Then,  $\mu^{\lambda_y}$ minimizes $\mu_1 \to D_{ni}(\mu_1\,||\, \mu_2)$, when $\mu_1 \in \mathcal{M}_\sigma$, and   $D_y^\sigma=- P(\log J_2 +\lambda_y f).$ Indeed, if $\int f d \mu_1=y$, we get that $ G (\mu_1)\leq G (\mu^{\lambda_y})=-P(\log J_2 + \lambda_y f).$  From Section 4.7 in \cite{Lop} we get  \eqref{triuvi}; it is a form of dynamical tilting probability.  \end{proof}
Now, consider the convex set 
$\Pi_{y}= \{\mu_1 \in \mathcal{M}\,|\, \int f d \mu_1=y \}$.

Given $y$, for $\mu_2$ fixed,
denote $D_y= D_{ni} (\Pi_y, \mu_2)= \inf \{D_{ni}(\mu_1, \mu_2)\,|\,\mu_1\in \Pi_y\}.$ 	A probability  $\mu_1$ minimizing $ D_y$ will be called the ni-Information projection of the triple $\mu_2$, $f$ and $y$. 
The next lemma considers a form of Canonical Ensemble probability (or, tilting probability in the terminology of \cite{Kosy}); we compute the expected value function and  entropy.

\begin{lemma} \label{ms1}  Given  $J_2 \in \mathfrak{J}$, consider the associated Gibbs state $\mu_2$ in $\mathcal{G}$, and a continuous function $B:\Omega \to \mathbb{R}.$  Consider a probability $\nu\in \mathcal{M}$  (not $\sigma$-invariant), in such a way that for any continuous function $g:\Omega \to \mathbb{R}$
	\begin{equation} \label{keq}\int g d \nu=\frac{1}{\int  e^{B \circ \sigma }d \mu_2}\,\int g   e^{B \circ \sigma }d \mu_2.	\end{equation}
Then,  for each $j$,
	$ E_{\nu} ( I_{\overline{j}} | \mathcal{X}_1)(x_1,x_2,..)= $
	\begin{equation} \label{jiu}  (\mathcal{L}_{J_2} (I_{\overline{j}}))\sigma(x)=J_2(j,x_2,x_3,..)= E_{\mu_2} ( I_{\overline{j}} | \mathcal{X}_1)(x_1,x_2,..).
	\end{equation}
	\end{lemma}
	
	\begin{proof}
	For any continuous function $w :\Omega \to \mathbb{R}$, as $\mu_2\in \mathcal{M}_\sigma$,

	$$\int w \circ \sigma I_{\overline{j}}d \nu = \frac{1}{\int  e^{B}d \mu_2} \int w \circ \sigma e^{B \circ \sigma} I_{\overline{j}}d \mu_2 = \frac{1}{\int  e^{B}d \mu_2} \int w e^B\, (\mathcal{L}_{J_2} \,(I_{\overline{j}})) d \mu_2=$$
	$$\frac{1}{\int  e^{B}d \mu_2} \int (w \circ \sigma)\, (\mathcal{L}_{J_2} \,(I_{\overline{j}}) \circ \sigma )e^{B \circ \sigma} d \mu_2=\int ( w \circ \sigma)\, (\mathcal{L}_{J_2} \,(I_{\overline{j}}) \circ \sigma\,) d \nu.$$
	
	Then, $J_2(j,x_2,x_3,..)= E_{\nu} ( I_{\overline{j}} | \mathcal{X}_1)(x_1,x_2,..)$. Moreover, as $\mu_2\in \mathcal{M}_\sigma$:
	
	$$h_{ni}(\nu) =- \sum_{j=1}^d \int  (\mathcal{L}_{J_2} \,(I_{\overline{j}}) \circ \sigma )\log (\mathcal{L}_{J_2} \,(I_{\overline{j}}) \circ \sigma )\, d \nu=$$
	\begin{equation} \label{kreq}-\int \sum_{j=1}^d   \mathcal{L}_{J_2} \,(I_{\overline{j}})  \,\log (\mathcal{L}_{J_2} \,(I_{\overline{j}})  )\,\frac{e^{B}}{\int e^{B  }d \mu_2} d \mu_2=\end{equation}
	$$	  - \sum_{j=1}^d \int J_2 (j, x_2,x_3,..)\log  J_2 (j, x_2,x_3,...)\, d \nu (x_1,x_2,..).$$
	The value \eqref{kreq} should be compared with the entropy of $\mu_2$, which is
	$$-\int \sum_{j=1}^d   \mathcal{L}_{J_2} \,(I_{\overline{j}})  \,\log (\mathcal{L}_{J_2} \,(I_{\overline{j}})  )\, d \mu_2.$$

	Taking into account \eqref{kreq}, the maximal (respectively, minimal) value of
	$$-\sum_{j=1}^d   \mathcal{L}_{J_2} \,(I_{\overline{j}})  \,\log (\mathcal{L}_{J_2} \,(I_{\overline{j}})  ),$$
	one can choose $B$ to enhance (respectively diminish)  $h_{ni}(\nu)$ with respect to  $h_{ni}(\mu_2)$.  This probability  $\nu$  may have larger or smaller ni-entropy than $\mu_2$.
	
From \eqref{jiu} we get an example of $\nu\ll\mu_2$ satisfying
$  \,\, D_{ni}(\mu_2\,||\, \nu) =$
	\begin{equation} \label{jiu1er}\int \sum_{j=1}^d I_{\overline{j}} \log ( E_{\mu_2} ( I_{\overline{j}} | \mathcal{X}_1)) d \mu_2  - \int \sum_{j=1}^d I_{\overline{j}} \log ( E_{\nu} ( I_{\overline{j}} | \mathcal{X}_1)) d \mu_2=0. \end{equation}
\end{proof}

Example \ref{orqi} generalizes Lemma \ref{ms1} by estimating the case $d\nu = \frac{ e^B}{\int e^B d\mu_2 } d\mu_2$.

\begin{theorem} Assume $\log J_2$ and $f:\Omega \to \mathbb{R}$ are H\"older functions. Given $y$, assume that $\lambda^*$ satisfies $\frac{d}{d\lambda}\log\int e^{\lambda f}\,d\mu_2=y$, where $\mu_2\in \mathcal{G}$ was obtained from the Jacobian $J_2$.
	Then,  an ni-Information Projection for  $\mu_2$, $f$, $y$, is
	\begin{equation} \label{triuvi87}d \nu_{\lambda^*}(x)= \frac{\sum_j\,e^{ \lambda^* f (j, \sigma (x))}J_2 (j, \sigma(x))  }{\int \sum_j\,e^{  \lambda^* f (j, \sigma (z))}J_2 (j, \sigma(z) ) d \mu_2(z)  } d \mu_2(x).  
	\end{equation}
	This probability maximizes $\mathcal{P}( \log J_2 + \lambda^* f).$

\end{theorem}

\begin{proof}
	$\lambda \to   \log(\int  e^{\lambda f} d\mu_2)$  has derivative $ \frac{d}{d \lambda}  \log(\int  e^{\lambda f} d\mu_2)|_{\lambda^*}=	\frac{\int f e^{\lambda^* \,f } d \mu_2 }{\int  e^{\lambda^* f} d \mu_2  }$.

	For $\mu_2$ fixed,  variables  $\lambda,k\in \mathbb{R}$, and a  probability $\mu_1$ in $\mathcal{M}$, consider
	$$ \sum_{j=1}^d\,  \,\log ( E_{\mu_2} (I_{\overline{j}}| \mathcal{X}_1) (x) )  E_{\mu_1} (I_{\overline{j}}| \mathcal{X}_1) (x) + \lambda\, \sum_{j=1}^d\,  \,  f(j, \sigma(x))E_{\mu_1} (I_{\overline{j}}| \mathcal{X}_1) (x) -  $$
	$$\sum_{j=1}^d\,   \log ( E_{\mu_1} (I_{\overline{j}}| \mathcal{X}_1) (x) )  \, E_{\mu_1} (I_{\overline{j}}| \mathcal{X}_1) (x) =
	$$
	$$ \sum_{j=1}^d\,  \,\log ( (J_2 \circ \sigma)(x) )  E_{\mu_1} (I_{\overline{j}}| \mathcal{X}_1) (x) + \lambda\, \sum_{j=1}^d\,  \,  f(j, \sigma(x))E_{\mu_1} (I_{\overline{j}}| \mathcal{X}_1) (x) +  $$
	\begin{equation} \label{klas7} - \sum_{j=1}^d\,   \log ( E_{\mu_1} (I_{\overline{j}}| \mathcal{X}_1) (x) )  \, E_{\mu_1} (I_{\overline{j}}| \mathcal{X}_1) (x) \end{equation}

	Integrating the above expression for $\mu_1\in \mathcal{M}$, we get the value
	\begin{equation}  \label{jiu1}- D_{ni} (\mu_1,\mu_2) + \lambda \int f d \mu_1.
	\end{equation}

To say $\mu_1$   maximizes \eqref{klas7} is equivalent to saying it maximizes
$\mathcal{P}( \log J_2 + \lambda f) 
	$

	Given $\lambda$, consider the $\nu_{\lambda}$  tilting  probability 
	\begin{equation} \label{tuvi}d \nu_{\lambda}(x)= \frac{\sum_j\,e^{ \lambda f (j, \sigma (x))}J_2 (j, \sigma(x))  }{\int \sum_r\,e^{  \lambda f (r, \sigma (z))}J_2 (r, \sigma(z) ) d \mu_2(z)  } d \mu_2(x).  \end{equation}
	
	From \eqref{jiu1er} we get $D_{ni}(\mu_2||\nu_{\lambda})=0.$

	 $\mu_2$ is $\sigma$-invariant and $\mathcal{L}_{J_2}^*(\mu_2)=\mu_2$.
	Suppose there exists a $\lambda^*=\lambda_{y}$ such that
	$$y= \int \frac{\sum_j\,f (j, \sigma (x))e^{ \lambda^* f (j, \sigma (x))}(J_2 (j, \sigma(x)) ) }{\int \sum_r\,e^{  \lambda^* f (r, \sigma (z))}(J_2 (r, \sigma(z)) d \mu_2(z)  } d \mu_2(x) =\frac{\int f e^{\lambda^* f} d\mu_2}{ \int \,e^{  \lambda^* f } d \mu_2 }.  $$

	Note that $\lambda^*$ can be obtained by solving $\frac{d}{d \lambda}  \log(\int f e^{\lambda f} d\mu_2)|_{\lambda^*}=y$ (differentiability follows  from \cite{PP}).
From Lemma \ref{ms1}, for any $j=1,2,...,d$,
	$$ E_{\nu_{\lambda^*}} (I_{\overline{j}}| \mathcal{X}_1) (x)= J_2(j,x_2,x_3,..) =J_2(j, \sigma(x)).$$

	For such $\lambda^*$, from \eqref{jiu1} and \eqref{jiu1er} we get $ - D_{ni} (\nu_{\lambda^*},\mu_2) + \lambda^* \int f d \nu_{\lambda^*} = \lambda^* \int f d \nu_{\lambda^*} =\lambda^*\, y .$ It is known that $D_{ni} (\mu_1,\mu_2)\geq 0$,  $\forall \mu_1\in \mathcal{M}$. In the case $\int f d \mu_1=y$, we get that $ - D_{ni} (\mu_1,\mu_2) + \lambda^* \int f d \mu_1= - D_{ni} (\mu_1,\mu_2) + \lambda^*  y \leq   - D_{ni} (\nu_{\lambda^*},\mu_2) + \lambda^* y.$
Therefore, $\nu_{\lambda^*}$ is  an ni-Information projection for  $\mu_2$, $f$, $y$. It is also a probability maximizing $\mathcal{P}( \log J_2 + \lambda^* f)$.
\end{proof}

\section{Results on the ni-pressure problem} \label{pre}
\begin{theorem} \label{PreA21q} 
Given the H\"older potential $A(x)= \log J(x)\leq 0$, $J\in \mathfrak{J}$,
we have $\mathcal{P}(A) = 0$. The supremum is attained by $\mu_J$, but is not unique in $\mathcal{M}$.
\end{theorem}
\begin{proof} 
	The probability $\mu_J$ satisfies $h_{ni}(\mu_J) = -\int \log J \, d\mu_J = -\int A \, d\mu_J$, hence $\mathcal{P}(A) \geq 0$. To show $\mathcal{P}(A) \leq 0$ for a generic $\rho \in \mathcal{M}$, we use the standard inequality $-\sum q_i \log q_i + \sum q_i \log p_i \leq 0$ (see \cite[Lemma 3.3]{PP}). 
	Setting $q_j = E_{\rho}(I_{\overline{j}} | \mathcal{X}_1)$ and $p_j = E_{\mu_J}(I_{\overline{j}} | \mathcal{X}_1) = J(j, x_2, \dots)$, we get:
	\begin{equation} \label{fgre1}
	-\sum_{j=1}^d E_{\rho}(I_{\overline{j}} | \mathcal{X}_1) \log E_{\rho}(I_{\overline{j}} | \mathcal{X}_1) + \sum_{j=1}^d E_{\rho}(I_{\overline{j}} | \mathcal{X}_1) \log J(j, x_2, \dots) \leq 0.
	\end{equation}
	Integrating with respect to $\rho$ directly yields $h_{ni}(\rho) + \int A \, d\rho \leq 0$. Example \ref{moj21} provides a non-invariant measure $\mu = \frac{1}{2}\delta_{1^\infty} + \frac{1}{2}\delta_{2, 1^\infty}$ where $h_{ni}(\mu)=\log 2$, demonstrating non-uniqueness.
\end{proof}
  
\begin{example} \label{foie} 
	In $\{1,2\}^\mathbb{N}$, let $\mu_1$ be a probability supported on $\overline{12}$. Its image $\rho = \mathcal{L}^*_{1/2}(\mu_1)$ is supported on $\overline{112}\cup \overline{212}$, and from Example \ref{moio53}, $h_{ni}(\rho)=\log 2$. For the potential $A=I_{\overline{112}\cup \overline{212}}$, the property $h_{ni}(\mu)\leq \log 2$ yields:
	\begin{equation} \label{fur3} 
		\mathcal{P}(A) = \sup_{\mu \in \mathcal{M}} \left( \int I_{\overline{112}\cup \overline{212}}(x) \, d\mu(x) + h_{ni}(\mu)\right) = 1 + \log 2.
	\end{equation}
	No invariant probability $\mu \in \mathcal{M}_\sigma$ attains this maximum. For $C=\overline{112}\cup\overline{212}$, and any $\mu \in 
\mathcal{M}_\sigma$ we get $\mu(C)=0$. This is so   because $C\cap\sigma(C)=\varnothing$, and from invariance $\mu(C)=\mu(\sigma^{-1}C)=0$, due to the fact that $\mu(C\cap\sigma(C)  )=0$.
 
 $C,\sigma^{-1}C,\sigma^{-2}C$ are pairwise disjoint, hence invariance gives $3\mu(C)\le1$. This  proves that no invariant measure can attain the value $1+\log2$. 
\end{example}

\begin{example} \label{foie1} 
	Consider the ni-pressure problem for $A = I_{\overline{1}} - I_{\overline{1}} \circ\sigma = I_{\overline{12}} - I_{\overline{21}}$ in $\{1,2\}^\mathbb{N}$. For any invariant $\mu \in \mathcal{M}_\sigma$, $\int A \, d\mu = 0$, so $\int A \, d\mu + h_{ni}(\mu) \leq \log 2$. 
	However, for the non-invariant probability $\mu = \frac{1}{2} \delta_{1,2,z} + \frac{1}{2} \delta_{2,2,z}$ (with fixed $z=(x_3,x_4,\dots)$), Example \ref{mojic} gives:
$ \int A \, d\mu + h_{ni}(\mu) = \frac{1}{2} + \log 2 > \log 2. $
	Thus, the maximizing probability for $\mathcal{P}(A)$ is strictly outside $\mathcal{M}_\sigma$.
\end{example}

\begin{example} \label{mai495}
	Let $J(1x)=p_1$ and $J(2x)=p_2$ where $p_1+p_2=1$ and $p_1 \neq 1/2$. Assume $p_1 \log p_1 < p_2 \log p_2$. As will be shown in Example \ref{mai49}, for $\mu_3 = \mathcal{L}_J^*(\mu_1)$, $h_{ni}(\mu_3) = -p_1 \log p_1 \mu_1(\overline{1}) - p_2 \log p_2 \mu_1(\overline{2})$. 
	For the potential $A=I_{\overline{1}}$, maximizing $\mathfrak{P}(I_{\overline{1}}) = \sup_{\mu \in \mathcal{M}} \left( \mu(\overline{1}) + h_{ni}(\mathcal{L}_J^*\mu) \right)$ is achieved simply by shifting all mass to $\overline{1}$. Thus, any $\mu$ with $\mu(\overline{2})=0$ attains the supremum.
\end{example}
Theorem \ref{tos} will show that maximizing $\mathfrak{P}(A)$ relates to maximizing a constrained dynamical capacity $H(\mathcal{L}_J^*, A, k)$. 

\section{Minimal and Maximal Output Entropy} \label{sec:Chan Cap}
Item 1. in the next Definition \ref{defimax31} is precisely expression (4.4) in \cite{King} for our setting.
From another point of view,  it is a distillation of an analogous concept considered in expression (2.9) in \cite{Faist} for the case of a quantum channel acting on density  matrices, and this is in consonance with our framework.

In the context of Blackwell's comparison of experiments \cite{Bla}, a family  (or experiment) $\mu_\theta^2$ is considered less informative than $\mu^1_\theta$ if there exists a channel $\Phi$ such that $\Phi(\mu^1_\theta) = \mu^2_\theta$; here families are indexed in $\mathcal{M}$ or $\mathcal{M}_\sigma$.

\begin{definition} \label{defimax31}\label{defimin}\label{defimaxd}\label{demin}\label{defimaxx}
	For a channel $\Phi$, we define the following Minimal and Maximal Output Entropy  problems:
	\begin{itemize}
		\item \textbf{1. Over all measures $\mathcal{M}$:} Minimal $H_\mathcal{M}(\Phi) = \inf_{\rho\in \mathcal{M}} h_{ni}(\Phi (\rho))$ and Maximal $H^\mathcal{M}(\Phi) = \sup_{\rho\in \mathcal{M}} h_{ni}(\Phi (\rho))$.
		\item \textbf{2. Over invariant measures $\mathcal{M}_\sigma$:} Dynamical Minimal $H_{\sigma}(\Phi) = \inf_{\rho\in \mathcal{M}_\sigma} h_{ni}(\Phi (\rho))$ and Dynamical Maximal $H^{\sigma}(\Phi) = \sup_{\rho\in \mathcal{M}_\sigma} h_{ni}(\Phi (\rho))$.
		\item \textbf{3. Over Gibbs measures $\mathcal{G}$:} Gibbs Minimal $H_\mathcal{G}(\Phi) = \inf_{\rho\in \mathcal{G}} h_{ni}(\Phi (\rho))$ and Gibbs Maximal $H^{\mathcal{G}}(\Phi) = \sup_{\rho\in \mathcal{G}} h_{ni}(\Phi (\rho))$.
	\end{itemize}
\end{definition}

$ H_\mathcal{M}(\mathcal{L}_J^*$) is called the minimal exit entropy for $\mathcal{L}_J^*$; it shows how noisy the channel is.
A probability $\tilde{\mu}\in \mathcal{M}$ minimizing $ H_\mathcal{M}(\mathcal{L}_J^*)$ is called the minimal input probability for the channel.
In Classical Information Theory, given an input probability $\rho$, the entropy of the output probability $\Phi (\rho)$
measures the uncertainty in what comes out. Low output entropy means the output is close to deterministic. A channel with small minimal output entropy preserves distinguishability well. The maximum output entropy of a channel  
tells you how spread out an output distribution can become. For some $J\in \mathfrak{J}$,  $\log d$ may not be realized for  the action  $\Phi= \mathcal{L}_J^*$ (Example \ref{mai49}).

\begin{remark} \label{ate} 
	By \eqref{moj47}, evaluating $H_{\sigma}(\mathcal{L}_J^*)$ is equivalent to finding minimizing probabilities on $\mathcal{M}_\sigma$ for the potential $B(x) = - \sum_{j=1}^d J(j,x)\log J(j,x)$, framing it strictly as an Ergodic Optimization problem.
\end{remark}

 Problem $H_\sigma(\mathcal{L}_J^*)$ effectively computes a \textit{dynamical} capacity of the channel $\mathcal{L}_{J}^*$ (over $\mathcal{M}_\sigma$), corresponding to scenarios where mutual information is zero. Example \ref{mai49} exhibits a case where the unique measure maximizing $H_\sigma(\mathcal{L}_J^*)$ in $\mathcal{M}_\sigma$ is \textit{not} $\mu_J$, but rather a measure supported on a fixed point for $\sigma$. 

\begin{proposition} \label{glq} 
	$H_\mathcal{G}(\Phi)=H_\sigma(\Phi)$.
\end{proposition}
\begin{proof} 
	For $J\in\mathfrak{J}$, the map $\mathcal{L}_{J}^* : \mathcal{M} \to \mathcal{M}$ is continuous. Since $h_{ni}$ is upper semicontinuous, the composition $\mu \to h_{ni}(\mathcal{L}_{J}^*(\mu))$ is also upper semicontinuous. Any $\rho\in\mathcal{M}_\sigma$ can be weakly approximated by elements in $\mathcal{G}$ (see \cite{Lo2}). Suppose $\rho$ attains the minimum of $H_\sigma$, and $\mu_n\in \mathcal{G}$ converges to $\rho$ as $n \to \infty$. Then, $\limsup_{n \to \infty} h_{ni}(\mathcal{L}_{J}^*(\mu_n)) \leq h_{ni}(\mathcal{L}_{J}^*(\rho))$, and the claim follows.
\end{proof}

\begin{proposition} \label{pqa}
	For $\Phi = \mathcal{L}_J^*$, Dirac masses realize extrema on $\mathcal{M}$:
	\begin{align}
		H_\mathcal{M}(\mathcal{L}_J^*) &= \inf_{x \in \Omega} \left( - \sum_{j=1}^d J(j,x)\log J (j, x) \right), \label{plus1} \\
		H^\mathcal{M}(\mathcal{L}_J^*) &= \sup_{x \in \Omega} \left( - \sum_{j=1}^d J(j,x)\log J (j, x) \right). \label{plus2}
	\end{align}
\end{proposition}
\begin{proof} 
	By \eqref{moj47}, $\mu_1 \to h_{ni}(\mathcal{L}_J^*(\mu_1))$ is affine in $\mu_1$. Therefore, the extrema of \eqref{plus2} are attained at $\delta_x$,  where $x$ is a maximum of $- \sum_{j=1}^d J(j,x)\log J (j, x)$. Problem \eqref{plus1} can be handled in a similar way.
\end{proof}

For a H\"older potential $A$, we can define constrained output entropy problems over $\mathfrak{K}(k)= \{\mu \in \mathcal{M} \mid \int A \, d\mu = k\}$:
\begin{align}
	H_s(\mathcal{L}_J^*, A, k) &= \sup_{\rho\in \mathfrak{K}(k)} h_{ni}(\mathcal{L}_J^*(\rho)), \label{pluu2} \\
	H_i(\mathcal{L}_J^*, A, k) &= \inf_{\rho\in \mathfrak{K}(k)} h_{ni}(\mathcal{L}_J^*(\rho)). \label{pluuu2}
\end{align}

For the constrained supremum $H_s$ the existence of a maximizer follows from compactness of the level set $\mathfrak K(k)$ and upper semicontinuity of $h_{ni}\circ\mathcal L_J^*$. 

\begin{theorem}\label{tos} 
	Given the ni-pressure problem $\mathfrak{P}(A)$, there exists $k$ such that the probability maximizing $H_s(\mathcal{L}_J^*,A,k)$ also maximizes $\mathfrak{P}(A)$.
\end{theorem}
\begin{proof} 
	Let $\tilde{\mu}_A \in \mathcal{M}$ be the maximizer of $\mathfrak{P}(A)$ and set $k = \int A \, d\tilde{\mu}_A$. If $\nu \in \mathfrak{K}(k)$ attains the supremum for $H(\mathcal{L}_J^*, A, k)$, then:
	$$ \mathfrak{P}(A) = \int A \, d\tilde{\mu}_A + h_{ni}(\mathcal{L}_J^*(\tilde{\mu}_A)) \leq \int A \, d\nu + h_{ni}(\mathcal{L}_J^*(\nu)). $$
As $\tilde{\mu}_A$ is a global maximizer, equality holds, and $\nu$ also maximizes $\mathfrak{P}(A)$.
\end{proof}

\begin{proposition}\label{toso} 
	Similarly, for the ni-pressure $\mathfrak{P}(-A)$, there exists $k\in \mathbb{R}$, such that the probability minimizing $H_i(\mathcal{L}_J^*,A,k)$ also maximizes $\mathfrak{P}(-A)$.
\end{proposition}

\begin{definition} \label{defica6} \label{defica3xz}
	For $\Phi = \mathcal{L}_J^*$, we define two concepts of  capacity:
	\begin{itemize}
		\item \textbf{Channel Information Capacity (CIC)} over $\mathcal{M}$:
		\begin{equation} \label{po01zz} 
			CIC_{\mathcal{M}} = \sup_{\mu \in \mathcal{M}} \left\{ h_{ni}(\Phi(\mu)) - \int h_{ni}(\Phi(\delta_x)) \, d\mu(x) \right\}.
		\end{equation}
		\item \textbf{Dynamical Channel Information Capacity (DCIC)} over $\mathcal{M}_\sigma$:
		\begin{equation} \label{po01m}
			DCIC_{\mathcal{M}_\sigma} = \sup_{\mu \in \mathcal{M}_\sigma} \left\{ h_{ni}(\Phi(\mu)) - \int h_{ni}(\Phi(\delta_x)) \, d\mu(x) \right\}.
		\end{equation}
	\end{itemize}
\end{definition}

The Channel Information Capacity  addresses the issue:
what is the maximum rate (in bits per channel use) at which information can be transmitted through a noisy channel with vanishing probability of error?  Here $\log$ denotes the natural logarithm, so $h_{ni}$ is measured in nats rather than in bits, as is common in Information Theory; \eqref{po01zz} and \eqref{po01m} are the analogues.

In Subsection \ref{anot}, using a different point of view,  we define another concept of  capacity  via \eqref{oca2},  presenting examples where the capacity is not zero.

\begin{proposition} \label{efica3}
	For $\Phi = \mathcal{L}_J^*$,  $DCIC_{\mathcal{M}_\sigma}$ and  $CIC_{\mathcal{M}}$ are both zero.
\end{proposition}

\begin{proof}Let $b(x)=-\sum_jJ(j,x)\log J(j,x)$. Theorem \ref{rreT} gives $h_{ni}(\Phi(\mu))=\int b\,d\mu$, and directly $h_{ni}(\Phi(\delta_x))=b(x)$. Hence $\int h_{ni}(\Phi(\delta_x))\,d\mu(x)=\int b\,d\mu=h_{ni}(\Phi(\mu))$, $\forall\mu\in\mathcal M$. Then, \eqref{po01zz} and \eqref{po01m} imply $CIC_{\mathcal M}=DCIC_{\mathcal M_\sigma}=0$.
\end{proof}

\begin{definition} \label{defica3vu} \label{deficaoi3}
	For $\Phi = \mathcal{L}_J^*$ and $\mu_3 = \Phi(\mu_1)$, the \textbf{average Mutual Information} $MIIO(X:Y)(\mu_1)$ =
\begin{equation} \label{po01de}
		 \int \sum_{j=1}^d I_{\overline{j}} \mathcal{L}_J(\log J) \, d\mu_1 - \int \sum_{j=1}^d I_{\overline{j}} \log(E_{\mu_3}(I_{\overline{j}} \mid \mathcal{X}_1)) \, d\mu_3.
	\end{equation}
	
\end{definition}
	\begin{proposition} For any $\mu_1\in \mathcal{M}$, we get $MIIO(X:Y)(\mu_1)=0.$ 
	
	\end{proposition}
	\begin{proof} With \eqref{po01de}, the first term is $\int\mathcal L_J(\log J)\,d\mu_1$, while the second term is $-h_{ni}(\Phi(\mu_1))$. Theorem \ref{rreT} makes these equal, so $MIIO(X:Y)(\mu_1)=0$.\end{proof}

\begin{remark} \label{tyr}

Denoting entropy by $S$, and  $\Phi=\mathcal{L}_J^*$ acting on $\mu\in \mathcal{M}$, {\it the Entropy Upper Bound for the channel} $\Phi$ (see Section  1 in  \cite{Amosov}) is the value
$$ EUB(\Phi)  = \sup_{\sum_{j=1}^r \pi_j=1, \pi_j\ge0, \mu_j\in \mathcal{M},\,r \in \mathbb{N} } \{S(  \sum_{j=1}^r \pi_j\, \Phi(\mu_j))-    \sum_{j=1}^r \pi_j\, S(\Phi(\mu_j)) \}.$$

	The EUB for a channel measures the convexity defect of the entropy. Because $h_{ni}(\mathcal{L}_J^*(\mu))$ is affine in $\mu$ (from \eqref{moj47}), the EUB for $\Phi = \mathcal{L}_J^*$ is  zero.
\end{remark}

\section{Applications} \label{sec:Examples}

This section provides several applications illustrating exotic phenomena within our  framework. First we compute the action of $\mathcal{L}_{J}^*$ on a general $\mu_1  \in \mathcal{M}$.

\begin{example} \label{mai9x}
	Let $\mu_1$ be the equilibrium probability for a Jacobian $J_1$, and  $\mu_3 = \mathcal{L}_{J}^*(\mu_1)$, which is generally not $\sigma$-invariant. For any cylinder $\overline{a_1, \dots, a_n}$:
	\begin{equation} \label{moj0}
		\mu_3(\overline{a_1, \dots, a_n}) = \int_{\overline{a_2, \dots, a_n}} J(a_1, y) \, d\mu_1(y).
	\end{equation}
	If $J = 1/d$ (corresponding to the maximal entropy measure $\mu$), equation \eqref{moj0} implies $\mu_1(\overline{a_2, \dots, a_n}) = (1/d)^{n-1}$. Thus, $\mathcal{L}_{J}^*(\mu_1) = \mu$ only if $\mu_1 = \mu$.
\end{example}

\subsection{The ni-Jensen-Shannon divergence} \label{JS}
We will introduce the
Jensen-Shannon divergence $JS(\mu_1,\mu_2)$, for $\mu_1,\mu_2\in \mathcal{M},$ and we will demonstrate that it is a pseudometric. We will show that it is a metric when restricted to $\mathcal{G}.$ \cite{Li} states that the concept of divergence introduced therein does not require assuming absolute continuity of the probabilities; rather, they are densities with respect to the counting measure on a finite set. Thus, it does not account for our framework.
Our definition in fact computes the distance between the expected value functions $E_{\,\mu\,} ( I_{\overline{j}} | \mathcal{X}_1)$ (and not between probabilities).

For simplicity, assume $\Omega=\{1,2\}^\mathbb{N}.$
For $0\leq p\leq 1$, $0\leq q\leq 1$, we set 
\begin{equation} \label{lkr} V(p,q) = p \log (\frac{2 p}{ p+q} ) + q \log (\frac{2 q}{ p+q} )\geq 0.\,\text{Then,}\, V(p,q)=0 \Rightarrow \,p=q. \end{equation}

Consider the probabilities $a=(a_1,a_2)$ and $b=(b_1,b_2)$, and denote
$$\alpha(a,b)=\{a_1 \log[a_1] + 
 a_2 \log[a_2] - ( a_1 \log[ 1/2  (a_1 + b_1)] + a_2 \log[ 1/2 (a_2 + b_2)])+$$
$$b_1 \log[b_1] + 
 b_2 \log[b_2] - ( b_1 \log[ 1/2  (a_1 + b_1)] + b_2 \log[ 1/2 (a_2 + b_2)])\}^{1/2}.$$
 
 From \eqref{lkr},  
 \begin{equation} \label{klm}\text{if}\,\, \alpha(a,b)=0,\,\, \text{then}\,\, a=b.
 \end{equation}
 
\begin{lemma} \label{potu}$\alpha$ satisfies the triangle inequality
\begin{equation} \label{kle} \alpha(a,c)\leq \alpha(a,b) + \alpha(b,c).
\end{equation}
\end{lemma}
\begin{proof} For  fixed $a,c$,  consider the function
\begin{equation} \label{klo}b_1 \to \beta(b_1)= \alpha(a,(b_1,1-b_1)) + \alpha((b_1,1-b_1),c)-  \alpha(a,c).\end{equation}
Using Mathematica one gets the explicit expression for the derivative $\frac{d}{d b_1}\beta(b_1)$ (a long expression).
W.l.o.g., we assume that $a_1<a_2$, $c_1<c_2$, and $a_1<c_1$.
At the points $b_1=a_1$ and $b_1=c_1$, there is a lack of differentiability, and  at these points we reach $\alpha(a,c)= \alpha(a,b) + \alpha(b,c)$. Moreover, there is a unique critical point: $\frac{d}{d b_1}\beta|_{(a_1+ c_1)/2}=0$, and it is a local maximum of $\beta$. The function $\beta$ is positive, except at the points $a_1$ and $c_1$; it is monotonically decreasing for $b_1<c_1$ and monotonically increasing for $c_1<b_1$.
\end{proof}

Given an {\it a priori probability} $\rho$ in $\Omega$, with support equal to $\Omega$, we set

$$  JS(\mu_1,\mu_2)=$$
$$ \int \{\,\sum_{j=1}^d  \, E_{\mu_1} ( I_{\overline{j}} | \mathcal{X}_1)[ \log ( E_{\mu_1} ( I_{\overline{j}} | \mathcal{X}_1))   -  \log \frac{1}{2}  (\, E_{\mu_1} ( I_{\overline{j}} | \mathcal{X}_1)) +   E_{\mu_2} ( I_{\overline{j}} | \mathcal{X}_1) \,) ]    + 
$$
$$\sum_{j=1}^d E_{\mu_2} ( I_{\overline{j}} | \mathcal{X}_1)[ \log ( E_{\mu_2} ( I_{\overline{j}} | \mathcal{X}_1))   -  \log \frac{1}{2}  (\, E_{\mu_1} ( I_{\overline{j}} | \mathcal{X}_1)) +   E_{\mu_2} ( I_{\overline{j}} | \mathcal{X}_1) \,) ] \}^{1/2} d \rho. $$

$JS$ is symmetric, and for general probabilities $\mu_1,\mu_2\in \mathcal{M}$ we get from \eqref{lkr} that $K_\rho(\mu_1,\mu_2)=0$ implies for all $j$ that  $E_{\mu_1} ( I_{\overline{j}} | \mathcal{X}_1)(x)= E_{\mu_2} ( I_{\overline{j}} | \mathcal{X}_1)(x)$, where they are positively defined.   In the case $\mu_1,\mu_2\in \mathcal{G}$, this  implies $\mu_1=\mu_2$. Given $\nu,\mu_2$  one can get  $K_\rho(\mu_2,\nu)=0$, but $\mu_2\neq \nu$ (see Lemma \ref{ms1}). The triangle inequality for $JS$ follows from Lemma \ref{potu}. Therefore, $K_\rho$ is a metric on $\mathcal{G}.$
When choosing  $\rho$ as the maximum entropy measure, this can be understood from a Bayesian perspective: given a lack of prior knowledge, the prior assumes that, in principle, there is no preference between $\mu_1,\mu_2$.


\subsection{The ni-Kullback-Leibler divergence and a  global non-perturbative bound}

Assume $\mu_1\in \mathcal{M}$ and $\mu_2 \in \mathcal{G}$, with Jacobian $J_2$.
Consider a H\"older function $f:\Omega \to \mathbb{R}$  and
take $\mu_\theta$, $\theta \in [0, \infty)$,  as the probability $d \mu_\theta= \frac{e^{\theta f}}{ \int e^{\theta f} d \mu_2} \, d \mu_2.$ 

Here it is natural to consider the cumulant generating function for $\mu_2$ as 
\begin{equation} \label{cum}\lambda(\theta)=\int  \log (\sum_{k} J_2(k x_2 x_3...) \, e^{\theta f(k x_2 x_3 ...)}\,) \,d \mu_1(x):=\log Z(\theta).
\end{equation}

We assume  that $\log ( E_{\mu_1} (I_{\overline{j}}| \mathcal{X}_1))$ is $\mu_1$-integrable, and we  will prove that
$\hspace*{4cm}    D_{ni}(\mu_1\,||\, \mu_2)- D_{ni}(\mu_1 \,||\, \mu_\theta) =$
\begin{equation} \label{cum1} \theta\, E_{\mu_1}(f)-  \int  \log (\sum_{k} J_2(k x_2 x_3...) \, e^{\theta f(k x_2 x_3 ...)}\,) \,d \mu_1(x)=\theta\, E_{\mu_1}(f)  - \lambda(\theta),\end{equation}
which is zero when $\theta=0.$
We also get a global non-perturbative upper bound
$$D_{ni}(\mu_1\,||\, \mu_2)- D_{ni}(\mu_1 \,||\, \mu_\theta) \leq$$
\begin{equation} \label{cum2} \theta [E_{\mu_1}  (f)  -  \int    \sum_{k} J_2(k x_2 x_3...) \, f(k x_2 x_3 ...))\, \,d \mu_1(x)].\end{equation}

From Jensen's inequality, \eqref{lort}, and \eqref{rw}, we get
$\,\,D_{ni}(\mu_1\,||\, \mu_2)- D_{ni}(\mu_1 \,||\, \mu_\theta) =$
$$  \sum_{j=1}^d(\int  E_{\mu_1} (I_{\overline{j}}| \mathcal{X}_1)\log ( E_{\mu_1} (I_{\overline{j}}| \mathcal{X}_1)) d \mu_1 (x)-\int \log J_2( \sigma(x)) d \mu_1(x))-$$
$$  \sum_{j=1}^d(\int  E_{\mu_1} (I_{\overline{j}}| \mathcal{X}_1)\log ( E_{\mu_1} (I_{\overline{j}}| \mathcal{X}_1)) d \mu_1 (x)-$$
$$\int    \sum_{j=1}^d  I_{\overline{j}}	(x)\log  (\, \frac{J_2(j x_2 x_3...) \, e^{\theta f (j x_2 x_3 ....)}}{\sum_{k} J_2(k x_2 x_3...) \, e^{\theta f(k x_2 x_3 ...)}}\,) d \mu_1(x))=$$
$$ \int    \sum_{j=1}^d  I_{\overline{j}}	(x)\log  (\, \frac{J_2(j x_2 x_3...) \, e^{\theta f (j x_2 x_3 ....)}}{\sum_{k} J_2(k x_2 x_3...) \, e^{\theta f(k x_2 x_3 ...)}}\,) d \mu_1(x) -$$
$$\int\sum_{j=1}^d  I_{\overline{j}}(x) \log J_2 (j x_2 x_3... )d \mu_1(x))=$$
$$ \int    \sum_{j=1}^d  I_{\overline{j}}	(x)\log  (\, \frac{\, e^{\theta f (j x_2 x_3 ....)}}{\sum_{k} J_2(k x_2 x_3...) \, e^{\theta f(k x_2 x_3 ...)}}\,) d \mu_1(x)=$$
$$ \sum_{j=1}^d \int      I_{\overline{j}}	(x)[ \log(e^{\theta f (j x_2 x_3 ....)}) - \log (\sum_{k} J_2(k x_2 x_3...) \, e^{\theta f(k x_2 x_3 ...)}\,) ]\,d \mu_1(x)=$$
$$\theta \int f d \mu_1  - \int   \log (\sum_{k} J_2(k x_2 x_3...) \, e^{\theta f(k x_2 x_3 ...)}\,) \,d \mu_1(x)\leq$$	 
$$\theta \int f d \mu_1  - \int    \sum_{k} J_2(k x_2 x_3...) \,\log( e^{\theta f(k x_2 x_3 ...)})\, \,d \mu_1(x)  =$$
$$\theta \int f d \mu_1  - \theta \int    \sum_{k} J_2(k x_2 x_3...) \, f(k x_2 x_3 ...))\, \,d \mu_1(x) =$$
$$\theta [\int f d \mu_1  -  \int    \sum_{k} J_2(k x_2 x_3...) \, f(k x_2 x_3 ...))\, \,d \mu_1(x)].$$

The function $\theta \to \lambda(\theta)$ is differentiable (see \cite{PP}),
\begin{equation} \label{cum3}\frac{d}{d \theta} \lambda(\theta)=\int   \frac{\sum_{k} f(k x_2 x_3 ...)\, J_2(k x_2 x_3...) \, e^{\theta f(k x_2 x_3 ...)}\,)}{ \sum_{k} J_2(k x_2 x_3...) \, e^{\theta f(k x_2 x_3 ...)} } \,d \mu_1(x) 
\end{equation}
and
$$\frac{d}{d \theta} \lambda(\theta)|_{\theta=0}=\int   \sum_{k} f(k x_2 x_3 ...)\, J_2(k x_2 x_3...) \, \,d \mu_1(x). $$

The Fisher information is
$I(\theta)  =\frac{d^2}{d \theta^2} \lambda(\theta).$

\subsection{Another point of view: mutual information and capacity} \label{anot}

In this subsection we follow a different point of view: the channel $\mathcal{L}_J^* $ will act on probabilities $q\in \Lambda$ and not on probabilities in $\mathcal{M}.$
Assume  $J\in \mathfrak{J}$ and $\mu \in \mathcal{G}$ is the unique one such that $\mathcal{L}_J^* (\mu)=\mu$.
By abuse of language, consider the extension of $q$ in $\Lambda$  to $\Omega$  via $q(y_1,y_2,...,y_n,...)=q(y_1).$
We adapt the notation of Section 12.4 in \cite{Moser}, denoting, for fixed $y\in \Omega$, $p_y(x)=W(y|x)=J(x*y)$, $x\in\{1,2,...,d\}$; this defines a probability $p_y\in \Lambda$.
In line with (12.45) in \cite{Moser}, 
the mutual information for the channel $\mathcal{L}_J^*$ acting  on a probability $q\in \Lambda$ is
$$I(X|Y):= I(q|W)=  \int \sum_{x=1}^d q(x) W(y|x) \log ( 
\frac{W(y|x) }{ \sum_{x^\prime} q(x^\prime) W(y|x^\prime)})d \mu(y)=$$
\begin{equation}\label{oca1} \int \sum_{x=1}^d q(x) J(x*y) \log ( \frac{J(x*y) }{ \sum_{x^\prime} q(x^\prime) J(x^\prime*y)})d \mu(y).
\end{equation}

The capacity of the channel $\mathcal{L}_J^*$ is by definition
\begin{equation} \label{oca2} 
C=\sup_{q\in \Lambda} I(q|W).
\end{equation}

A probability $\tilde{q}$ maximizing \eqref{oca2} is called an optimal input probability for the channel.

\begin{example} \label{oca4} Suppose $J=\frac{1}{d}$, so that $\mu$ is the measure of maximal entropy. Then, $\sum_{x^\prime} q(x^\prime) J(x^\prime*y)
= \frac{1}{d}$, independently of $q_0(x)>0$ or not. In this case $C=0$, and any probability $q \in \Lambda$ maximizes \eqref{oca2}. 
\end{example}

\begin{example} \label{oca4} Take $J=p_1 I_{\overline{1}} + p_2 I_{\overline{2}} $, and  $\mu$   the associated independent probability on $\Omega=\{1,2\}^\mathbb{N}$. Consider  the  function  $(q_1,q_2)\to I((q_1,q_2),W) $, using Mathematica  we get that the  maximum of $I$ is attained when $\tilde{q}=(\tilde{q}_1,\tilde{q}_2)$ is
$$\tilde{q}_1=\frac{e^{\frac{p_2-p_1 + p_1 \log(p_1) -p_2 \log (p_2)     }{p_1-p_2}}-p_2}{p_1-p_2}  \,\,\text{and} \,\,\tilde{q}_2=\frac{-e^{\frac{p_2-p_1 + p_1 \log(p_1) -p_2 \log (p_2)     }{p_1-p_2}}+ p_1}{p_1-p_2}.  $$

 In this way we get  an explicit expression  for $C=(\tilde{q},W).$ The maximal value of the capacity is $e^{-1}$ and is attained when $p_1=0$ or $p_1=1$, and the minimal value is zero and is attained when $p_1=1/2$. Moreover,  $e^{-1}\leq \tilde{q}_1\leq 1- e^{-1}.$

\end{example}

\subsection{Computing ni-entropy} \label{nini}

\begin{example} \label{mojic} 
	For $\Omega=\{1,2\}^\mathbb{N}$, let $\mu = p_1\delta_{1^\infty} + p_2\delta_{2, 1^\infty}$ with $p_1+p_2=1$. 
	Here, $E_{\mu}(I_{\overline{j}} | \mathcal{X}_1) = p_j  (I_{\{1^\infty\}} + I_{\{2 1^\infty\}})$, yielding:
	$ h_{ni}(\mu) = -p_1 \log p_1 - p_2 \log p_2 > 0.$
	A similar property holds for any $\mu = p_1\delta_{x_1} + p_2\delta_{x_2}$ as long as $\sigma(x_1) = \sigma(x_2)$. If $p_1=1/2=p_2$, we get maximal entropy $ h_{ni}(\mu)=\log 2$.
\end{example}

\begin{example} \label{maili7}
	Conversely, consider $\mu = p_1\delta_{x_1} + p_2\delta_{x_2}$ where $x_1 = (1, r_2, \dots)$ and $x_2 = (2, s_2, \dots)$. If $\sigma(x_1) \neq \sigma(x_2)$, then $E_{\mu}(I_{\overline{1}} | \mathcal{X}_1) = I_{\{\sigma(x_2)\}} \circ \sigma $. One can take 	$ E_{\mu}(I_{\overline{2}} | \mathcal{X}_1) =I_{\{\sigma(x_1)\}}\circ \sigma $. It follows that $h_{ni}(\mu) = 0$. The dynamics of $\sigma$ fully discriminate whether the ni-entropy is positive (Example \ref{mojic}) or zero. This case covers the situation where $\mu\in \mathcal{M}_\sigma$ has support in a orbit $\{x_1=x\in \overline{1}, x_2=\sigma(x)\in \overline{2}\}$ of period $2$.  $E_{\mu}(I_{\overline{1}} | \mathcal{X}_1)(x_1)=0$, and as $x_1$  is in the support of $\mu$, it follows that $\log E_{\mu}(I_{\overline{1}} | \mathcal{X}_1)$ is not $\mu$-integrable.
\end{example}

\begin{remark} \label{hjq}
	\textbf{Discontinuity of expected values and ni-entropy.} Let $x \in \overline{1}$ and choose $x_k \in \overline{2}$ such that $\sigma(x) \neq \sigma(x_k)$ but $x_k \to y$, where $\sigma(x) = \sigma(y)$. For $\mu_k = p_1\delta_x + p_2\delta_{x_k}$, we have $h_{ni}(\mu_k) = 0$ (see Example \ref{maili7}). However, $\mu_k$ weakly converges to $\mu = p_1\delta_x + p_2\delta_y$, which has $h_{ni}(\mu) > 0$ (by Example \ref{mojic}). Thus, $h_{ni}$ is not  continuous on $\mathcal{M}$.
\end{remark}

\begin{example} \label{moic16} Take $\Omega=\{1,2,3\}^\mathbb{N}$, $\mu =\frac{1}{3} (\delta_{\{(112)^\infty\}} + \delta_{\{(121)^\infty\}} + \delta_{\{(211)^\infty\}}).$
Then, $E_{\mu}(I_{\overline{3}} | \mathcal{X}_1)=0 = 0 I_{\{(112)^\infty\}}+  0 I_{\{(121)^\infty\}} +0 I_{\{(211)^\infty\}},$ and $\log(E_{\mu}(I_{\overline{3}} | \mathcal{X}_1))$ is not $\mu$-integrable.
In this case $B=\sigma^{-1}(\overline{1})\in \mathcal{X}_1$ satisfies $I_{\overline{3}}=0$, $\mu$.a.e on $B$. Moreover, $E_{\mu}(I_{\overline{1}} | \mathcal{X}_1) =I_{ \overline{12} }\circ \sigma+ I_{ \overline{21} }\circ \sigma $, $E_{\mu}(I_{\overline{1}} | \mathcal{X}_1) =I_{ \overline{11} }\circ \sigma$ and $h_{ni}(\mu) = 0.$
\end{example}

\begin{example} \label{moj16} 
	For any Dirac measure $\mu = \delta_x$, we get $E_{\mu}(I_{\overline{j}} | \mathcal{X}_1) = I_{\{\sigma(x)\}}\circ \sigma$, if $x \in \overline{j}$, and zero otherwise. Consequently, $h_{ni}(\delta_x) = 0$ for all $x \in \Omega$.
\end{example}

\begin{example}
 From Example \ref{moj16}: $  E_{\rho } ( I_{\overline{1}}\, |\, \mathcal{X}_1)=I_{\{1^\infty\}} \circ \sigma:=r(x)$, when $\rho=\delta_{1^\infty}$.
For $\mu=\delta_{1 ,1 ,2,2,2,2,...}$, we get $  E_{\mu } ( I_{\overline{1}}\, |\, \mathcal{X}_1)=I_{\{1,2,2,2,...\}} \circ \sigma:= s(x)$, and $  E_{\mu } ( I_{\overline{2}}\, |\, \mathcal{X}_1)=0$.
Note that $\frac{1}{2}r +\frac{1}{2}s \neq  E_{\frac{1}{2}\rho +\frac{1}{2} \mu } ( I_{\overline{1}}\, |\, \mathcal{X}_1)=r+s $. When   $\rho$ and $\mu$ have disjoint supports  $ E_{\lambda\rho +(1-\lambda) \mu } ( I_{\overline{1}}\, |\, \mathcal{X}_1)= E_{\rho  } ( I_{\overline{1}}\, |\, \mathcal{X}_1)+  E_{\mu } ( I_{\overline{1}}\, |\, \mathcal{X}_1)$.
The analogous property does not hold for convex combinations of probability measures whose supports intersect each other.

\end{example}

\begin{example} \label{eigp}
	Let $\mu$ be an eigenprobability of the dual Ruelle operator $\mathcal{L}_B^*$ (i.e., $\mathcal{L}_B^*(\mu) = \lambda \mu$).  The  ni-entropy can be explicitly computed. Indeed,   by evaluating $\int (g \circ \sigma) I_{\overline{j}} \, d\mu$, $j=1,2,...,d$, one gets the  expected value function:
\begin{equation} \label{confa} E_{\mu}(I_{\overline{j}} | \mathcal{X}_1)(x) = \frac{e^{B(j, \sigma(x))}}{\sum_k e^{B(k, \sigma(x))}}. \end{equation}
\end{example}

\begin{example} \label{mai}
	For an independent process where $J(j, x) = p_j$, we simply get $E_{\mu}(I_{\overline{j}} | \mathcal{X}_1) = p_j$. Thus, $h_{ni}(\mu) = h(\mu) = -\sum p_j \log p_j$.
\end{example}

\begin{example} \label{orqi}
	Generalizing Lemma \ref{ms1}, we can compute $h_{ni}(\nu)$, where $d\nu =  \frac{e^B}{ \int e^B d \mu_2} d\mu_2$. The expected value function $W_j \circ \sigma = E_{\nu}(I_{\overline{j}} | \mathcal{X}_1)$ must satisfy:
	$$ W_j \, \mathcal{L}_{J_2}(e^B) = \mathcal{L}_{J_2}(I_{\overline{j}} \, e^B), $$
where $J_2$ is the Jacobian of $\mu_2$. 	Evaluating $\mathcal{L}_{J_2}$ explicitly yields: $\exists c>0$ s.t.
	\begin{equation} \label{rw}
	 E_{\nu}(I_{\overline{j}} | \mathcal{X}_1)=	(W_j \circ \sigma)(x) = \frac{J_2(j, \sigma(x)) \, e^{B(j, \sigma(x))}}{\sum_k J_2(k, \sigma(x)) \, e^{B(k, \sigma(x))}}>c>0.
	\end{equation}
\end{example}

The entropy of a non-stationary Markov Process.

\begin{example} \label{gost1} 
	Consider a $2 \times 2$ stochastic matrix $P$ (with all entries positive) and an initial probability vector $\pi =(\pi_1,\pi_2)$, defining a not necessarily stationary Markov Process, which produces a probability $\mu\in \mathcal{M}$ on $\Omega=\{1,2\} ^\mathbb{N}$. We can explicitly compute $h_{ni}(\mu)$ using the conditional expectations:
	\begin{equation}
		E_{\mu}(I_{\overline{1}} | \mathcal{X}_1)=c_1 (I_{\overline{11} } + I_{\overline{21} }) + c_2 ( I_{\overline{12} } + I_{\overline{22} } ),
	\end{equation}
	where $c_1= \frac{\pi_1 P_{11} }{\pi_1 P_{11} +\pi_2 P_{21} }$ and $c_2= \frac{\pi_1 P_{12} }{\pi_1 P_{12} +\pi_2 P_{22} }$. Similarly,
	\begin{equation}
		E_{\mu}(I_{\overline{2}} | \mathcal{X}_1)=d_1 (I_{\overline{11} } + I_{\overline{12} }) + d_2 ( I_{\overline{21} } + I_{\overline{22} } ),
	\end{equation}
	where $d_1= \frac{\pi_2 P_{21} }{\pi_1 P_{11} +\pi_2 P_{21} }$ and $d_2= \frac{\pi_2 P_{22} }{\pi_1 P_{12} +\pi_2 P_{22} }$.
\end{example}

\subsection{Computing ni-Kullback-Leibler divergence} \label{KLD}

\begin{remark} \label{noe}	  In Classical Information Theory, the strong data processing inequality ensures $D_{KL}(\Phi(p) \parallel \Phi(q)) \leq D_{KL}(p \parallel q)$ for distributions on finite sets. In our dynamical setting with $\Phi = \mathcal{L}_J^*$ and $\mu, \rho \in \mathcal{M}$, the analogous inequality $D_{ni}(\Phi(\mu) \parallel \Phi(\rho)) \leq D_{ni}(\mu \parallel \rho)$ does \textbf{not} always hold. 
\end{remark}

\begin{example} \label{reco}
	As in Example \ref{orqi}, let $\mu_2 \in \mathcal{G}$ be a Gibbs state for $J_2 \in \mathfrak{J}$, and define $\nu$ via $d\nu = \frac{ e^B}{\int e^B d \mu_2} d\mu_2$. We will show that $D_{ni}(\mu_2\,||\, \nu) = 0$. 
	
	The expected value $E_{\nu} ( I_{\overline{j}} | \mathcal{X}_1)$ is given by \eqref{rw}. Applying Jensen's inequality for convex functions to the $\mathcal{X}_1$-integrable function $\log ( E_{\nu} ( I_{\overline{j}} | \mathcal{X}_1))$, we get
$$-\log (\sum_{k} J_2(k ,x_2, x_3,...) \, e^{B(k, x_2, x_3 ...)})\,\leq$$
$$-\sum_{k} J_2(k, \sigma(x_1,x_2 x_3...) \, B(k,\sigma(x_1, x_2 x_3 ...)).$$
	
	Using the above inequality, one can show that
	$$ 0 \leq D_{ni}(\mu_2\,||\, \nu) = \int \sum_{j=1}^d I_{\overline{j}} \log ( E_{\mu_2} ( I_{\overline{j}} | \mathcal{X}_1)) d \mu_2 - \int \sum_{j=1}^d I_{\overline{j}} \log ( E_{\nu} ( I_{\overline{j}} | \mathcal{X}_1)) d \mu_2 \leq$$
	$$\int B(x)\,d\mu_2(x) - \int   \sum_{k=1}^d J_2(k, \sigma(x_1,x_2 x_3...) \, B(k,\sigma(x_1, x_2 x_3 ...))\, d \mu_2(x_1,x_2,...)\leq$$
	$$ \leq \int B(x)\,d\mu_2(x) - \int B(x)\,d\mu_2(x) = 0. $$
\end{example}

\begin{example} \label{orqi8}
	Given a Gibbs state $\mu_2$, consider two continuous functions $B_1 \neq B_2$ and define $d\nu_k= \frac{ e^{B_k}}{\int e^{B_k} d \mu_2} \mu_2$ for $k=1,2$. Then,
	\begin{equation}\label{rw7w}
		E_{\nu_k} ( I_{\overline{j}} | \mathcal{X}_1)=	(W_j^k\circ \sigma)(x) = \frac{J_2(j x_2 x_3...) \, e^{B_k(j x_2 x_3 ....)}}{\sum_{r=1}^d J_2(r x_2 x_3...) \, e^{B_k(r x_2 x_3 ...)}}.
	\end{equation}
	Considering $B_1 \neq B_2$,  one can show  that $D_{ni}(\nu_2\,||\, \nu_1) =$
	$$(\,\int B_1(x)\,d\nu_2(x) - \int \log (\sum_{r} J_2(r, x_2 ,x_3...) \, e^{B_1(r ,x_2 ,x_3, ...)}) \,d \nu_2(x_1,x_2,...)\,)\neq 0.$$
\end{example}

The above shows that it makes a difference to assume the existence of a dominating measure,  or to assume they are just absolutely continuous.

\begin{example} \label{eig_kl_example} Given the H\"older potential $B:\Omega=\{1,2\}^\mathbb{N} \to \mathbb{R}$ consider the corresponding Ruelle operator $\mathcal{L}_B$, and via the Ruelle Theorem, the associated eigenvalue $\lambda$ and eigenfunction $f$  (that is, $\mathcal{L}_B (f)=\lambda f$).
We denote by $\rho$ in $\mathcal{G}$ the equilibrium  probability for the classical pressure problem $P(B)$ (a certain $J \in \mathfrak{J}$ is the associated Jacobian). Then, $\mathcal{L}_J^*(\rho)=\rho$. Denote by $\mu$ the eigenprobability for $B$,  that is, $\mathcal{L}_B (\mu)^*=\lambda \mu$, as considered in Example \ref{eigp}. 

Using \eqref{confa}  and  the coboundary relation $e^{B(j, x)} = \frac{J(j,x) f(j,x)}{\lambda f(x)}$, one can get:
	$$ D_{ni} (\rho \, || \, \mu) = \int \log \left( \frac{ f(x)}{\lambda\, f(\sigma(x))\, \sum_{j=1}^d e^{B(j,\sigma x)}} \right) d \rho(x). $$

\end{example}

\begin{example}
\label{orqi45} 
Given $Q\in \Lambda, \mu\in \mathcal{G}$, set $\rho := Q \mu\in \mathcal{M}$, s. t.  $\forall$ continuous $f$
$$ \int f \,d \rho = \sum_{j=1}^d Q_j\int_{\overline{j}  }\frac{1}{\mu(\overline{j})}f \, d\mu= \int  \sum_{j=1}^d \frac{Q_j}{\mu(\overline{j})} I_{ \overline{j}}\,  f d\mu .$$
Then, for any $j=1,2,...,d$, one can show that
\begin{equation} \label{plo1} E_{\rho} ( I_{\overline{j}} | \mathcal{X}_1))(x)= \frac{ J (j,\sigma(x)) \frac{Q_j}{\mu(\overline{j})}}{\sum_{r=1}^d  J (r,\sigma(x)) \frac{Q_r}{\mu(\overline{r})} } .\end{equation}

  Given $Q^1, Q^2\in \Lambda$ with all entries positive, \eqref{plo1} for  $  Q_1 \mu$ and $Q_2 \mu$ is strictly positive, and the both supports are equal to $\Omega$. One can show that
 $$ D_{ni} (Q_1 \mu, Q_2 \mu)=0\,\, \Rightarrow \,\,   E_{Q_1 \mu} ( I_{\overline{j}} | \mathcal{X}_1))= E_{Q_2 \mu} ( I_{\overline{j}} | \mathcal{X}_1)).$$
 This does not mean that $Q_1=Q_2$.
 Now,  using \eqref{opayu}:
 $$0= \sum_j \int E_{Q_1 \mu} (I_{\overline{j}}| \mathcal{X}_1)\log ( E_{Q_1 \mu} (I_{\overline{j}}| \mathcal{X}_1))dQ_1 \mu- $$
 $$\sum_j\int E_{Q_1 \mu} (I_{\overline{j}}| \mathcal{X}_1)\log ( E_{Q_2\mu} (I_{\overline{j}}| \mathcal{X}_1))d Q_1 \mu=$$
 $$ \sum_a  Q^1_a \int \frac{1}{ \mu (\overline{a}) }   I_{\overline{a}}\sum_j   E_{Q_1 \mu} (I_{\overline{j}}| \mathcal{X}_1)\log ( E_{Q_1 \mu} (I_{\overline{j}}| \mathcal{X}_1))d \mu- $$
 $$   \sum_a  Q^1_a \int \frac{1}{ \mu (\overline{a}) }  I_{\overline{a}}  \sum_j  E_{Q_1 \mu} (I_{\overline{j}}| \mathcal{X}_1)\log ( E_{Q_2 \mu} (I_{\overline{j}}| \mathcal{X}_1))d \mu =$$
 $$ \sum_a  Q^1_a  \frac{1}{ \mu (\overline{a}) }  \int [\sum_j   E_{Q_1 \mu} (I_{\overline{j}}| \mathcal{X}_1)\log ( E_{Q_1 \mu} (I_{\overline{j}}| \mathcal{X}_1))- $$
 $$E_{Q_1 \mu} (I_{\overline{j}}| \mathcal{X}_1)\log ( E_{Q_2 \mu} (I_{\overline{j}}| \mathcal{X}_1))] (a, \sigma(x))d \mu(x)\geq 0.$$ 

Then, for any $a=1,2,...,d$, and $\mu$ a. e.  $x$, the continuous functions 
$$E_{Q_1 \mu} (I_{\overline{j}}| \mathcal{X}_1))(a \sigma(x))= E_{Q_2 \mu} (I_{\overline{j}}| \mathcal{X}_1))(a \sigma(x)) .$$

 \end{example}

\subsection{The action of the dual of the Ruelle operator} \label{runi}

$\mathcal{L}_J^*$ is not necessarily surjective and can either increase or decrease entropy.

\begin{example} \label{mai49}
	Let $J(1,x)=p_1$ and $J(2,x)=p_2$ where $p_1+p_2=1$ and $p_1 \neq 1/2$. Let $\mu$ be the associated independent probability. For any $\mu_1$, the image $\mu_3 = \mathcal{L}_J^*(\mu_1)$ maintains the expected values $E_{\mu_3}(I_{\overline{j}} | \mathcal{X}_1) = p_j$. Thus:
	\begin{equation} \label{moj473} 
		h_{ni}(\mu_3) = -p_1 \log (p_1) \mu_1(\overline{1}) - p_2\log (p_2) (1-\mu_1(\overline{1} )). 
	\end{equation}
	This value is strictly less than $\log 2$, meaning the maximal entropy measure is absent from the image of $\mathcal{L}_J^*$. If $\mu_1$ is an independent probability for $(q_1, q_2)$, $h_{ni}(\mu_3)$ can be strictly larger or smaller than $h_{ni}(\mu_1)$ depending on the choice of $q_1$, showing entropy can both increase and decrease. When considering maximization  of $h_{ni}(\mu_3)$ for $\mu_3\in\mathcal{M}$, many non-invariant probabilities supported on $\overline1$ or $\overline2$ attain the extrema.
For the maximization   of $h_{ni}(\mu_3)$ on $\mathcal{M}_\sigma$, 
the extremal values  are uniquely attained by the Dirac measures $\delta_{1^\infty}$ and $\delta_{2^\infty}$.

\end{example}

\begin{example} \label{oiaq} 
	On $\Omega=\{0,1\}^\mathbb{N}$, define a Jacobian $J_{i,j} = P_{ji}$ from a stochastic matrix $P$. For $\mu_1 = q_0 I_{\overline{0}}\mu_0 + q_1 I_{\overline{1}}\mu_0$ (where $\mu_0$ is the maximal entropy measure), evaluating limits on cylinder sets shows that for $\mu_3 = \mathcal{L}^*_J (\mu_1)$, the expected value is explicit: $ E_{\mu_3 } ( I_{\overline{0}}\, |\, \mathcal{X}_1) = (J_{00} I_{\overline{0}}+ J_{01} I_{\overline{1}}) \circ \sigma$;  $h_{ni}(\mu_3)$ follows.
\end{example}

\begin{example} \label{moio53} 
	When $J=1/d$, applying $\mathcal{L}^*_{1/d}(\mu_1)=\mu_3$, for which we get
	$$ h_{ni}(\mu_3) = - \int_{\Omega} \sum_{j=1}^d \frac{1}{d} \log \left(\frac{1}{d}\right) \, d\mu_1 = \log d. $$
	Thus, $\mathcal{L}^*_{1/d}$ never decreases entropy. Note that $\mu_3$ is not $\sigma$-invariant unless $\mu_1$ is already the maximal entropy measure $\mu$.
\end{example}

\begin{example} \label{mailr3}
	Let $J(j,x)=p_j$, $j=1,2,...,d$,  correspond to an independent probability $\mu$. If $\mu_1$ is an independent probability associated to $(q_1,q_2,...,q_d)$, then $\mu_3 = \mathcal{L}_{J}^* (\mu_1)$ is not $\sigma$-invariant (e.g., $\mu_3(\overline{a_1,\dots,a_n}) = p_{a_1} q_{a_2}\dots q_{a_n}$). However, $h_{ni}(\mu_3)$ evaluates to $-\sum p_j \log p_j = h_{ni}(\mu)$, showing $\mathcal{L}_{J}^*$ can map distinct measures $\mu_1$ into states with the same entropy.
\end{example}

\begin{example} \label{moj25} 
	Let $J=1/2$ and $\mu_1=\delta_{1^\infty}$. Then $\mu_3 = \mathcal{L}_{J}^* (\mu_1) = \frac{1}{2}\delta_{1^\infty}+ \frac{1}{2}\delta_{2, 1^\infty}$. From Example \ref{mojic}, $h_{ni}(\mu_3) = \log 2 > h_{ni}(\mu_1) = 0$. This confirms $\mathcal{L}_{J}^*$ can strictly increase entropy, mapping a Dirac measure to a non-invariant measure with positive entropy.
	Given $g\in \mathcal{C}$,   $\frac{1}{2}[ g(1^\infty) + g(2\, 1^\infty)]= \int g d \mu_3.$
	
\end{example}

\section{Conclusion} \label{Conc}

We introduced a natural definition for the entropy of probability measures on symbolic spaces that are not necessarily shift-invariant. This concept extends Kolmogorov-Sinai entropy to non-stationary settings and diverges from classical Shannon entropy by not requiring a fixed dominating measure on $\Omega=\{1,2,...,d\}^\mathbb{N}$. Consequently, it allows the thermodynamic pressure formalism to be broadened. Through various examples, we illustrated the complexity of this framework, particularly in analyzing the dual Ruelle operator's action as an information channel on the space of probability measures $\mathcal{M}$. We relate problems in Information Projections to the maximization of ni-pressure. We derived properties related to the concepts of ni-Jensen-Shannon divergence, mutual information and capacity. We enlarged the scope of Information Theory by not assuming the existence of a dominating probability. Future explorations of this framework hold significant potential, including applications of ni-divergence and ni-cross entropy in Bayesian Statistics,   Chernoff-type concentration bounds and Geometric Information. One can ask for applications of ni-divergence in areas like  maximum likelihood estimation and the Akaike Information Criterion in model selection. Given probabilities $\mu_2\in \mathcal{G}$ and $\nu$ as in Example \ref{orqi}, another possible  topic of investigation is   the ni-cross entropy,  naturally defined by
$$ H_{ni}(\mu_2,\nu) =  - \int \sum_{j=1}^d  E_{\nu} ( I_{\overline{j}} | \mathcal{X}_1)\log ( E_{\nu} ( I_{\overline{j}} | \mathcal{X}_1)) d \mu_2.
$$

In (2) in  \cite{Max} the author mentioned the existence of a constant $c$ such that
$$W_1 ( \nu, \mu_\phi)\leq \,c\, \sqrt{P(\phi)  - h(\nu) - \int \phi \,d \, \nu   },$$
where $W_1$ is the $1$-Wasserstein distance, $ P(\phi)$ is the pressure of $\phi$ (see \eqref{AA}), $\nu \in \mathcal{M}_ \sigma$, $\mu_\phi$ the classical equilibrium probability for $\phi$,   and $h(\nu)$ is   Kolmogorov entropy. We wonder if a similar expression may be true when $\nu\in \mathcal{M}$, replacing $h(\nu)$ by $h_{ni}(\nu),$ and
$ P(\phi)$ by $\mathcal{P} (\phi)$?

Statement of AI usage - The strategy followed in the paper to obtain the main results was carried out without the aid of AI. Claude was consulted - and its help obtained  - only for specific, isolated points concerning classical topics in measure theory and the topology of metric spaces, which are standard material found in textbooks on the subject. A very small number of examples (3 or 4 in total) were suggested by Claude
We take full responsibility for the accuracy and originality of the final submission.

\end{document}